\documentclass[11pt]{amsart}

\usepackage{amscd}
\usepackage{caption2}
\usepackage[dvips]{graphicx}
\usepackage{psfrag}
\usepackage[psamsfonts]{amssymb}
\usepackage{enumerate}
\usepackage{flafter}
\usepackage[all]{xy}

\newtheorem{Thm}{Theorem}[section]
\newtheorem{Lem}[Thm]{Lemma}

\newtheorem{defn}[Thm]{Definition}
\newtheorem*{claim}{Claim}

\theoremstyle{definition}

\newtheorem{Ques}[Thm]{Question}

\theoremstyle{remark}

\newtheorem{Exa}[Thm]{Example}

\newcommand{\Z}{{\mathbb Z}}

\newcommand{\R}{{\mathbb R}}

\newcommand{\norm}[1]{\left\Vert#1\right\Vert}
\newcommand{\abs}[1]{\left\vert#1\right\vert}
\newcommand{\set}[1]{\left\{#1\right\}}

\newcommand{\To}{\longrightarrow}

\newcommand{\dash}[1]{\stackrel{#1}{\mbox{---}}}
\def\co{\colon\thinspace}
\allowdisplaybreaks

\begin{document}

\title{Lectures on   topology of words}
\author[Vladimir Turaev]{Vladimir Turaev}
                     \address{%
              IRMA, Universit\'e Louis  Pasteur - C.N.R.S., \newline
\indent  7 rue Ren\'e Descartes \newline
                     \indent F-67084 Strasbourg \newline
                     \indent France \newline
\indent e-mail: turaev@math.u-strasbg.fr } \dedicatory{Based on
notes by Eri Hatakenaka, Daniel Moskovich, and  Tadayuki Watanabe}
\begin{abstract}  We discuss a topological approach to words introduced
  by the author in \cite{Tu1}--\cite{Tu3}.
  Words on an arbitrary alphabet are
   approximated by Gauss words and then studied
  up to natural modifications
inspired by the   Reidemeister moves on knot diagrams. This leads us to
a notion of homotopy for words. We  introduce several homotopy
invariants of  words and give a homotopy classification of  words of
length five.
                     \end{abstract}
             \maketitle

\section{Introduction}

\begin{figure}[htbp]
\psfrag{a}[cc][cc]{$a$}
\psfrag{b}[cc][cc]{$b$}
\psfrag{c}[cc][cc]{$c$}
\psfrag{d}[cc][cc]{origin}
\includegraphics[scale=0.6]{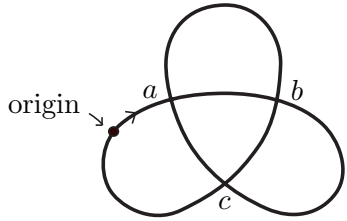}
\caption{A plane curve associated with the word $abcabc$}
\label{fig1}
\end{figure}

Words are finite sequences of symbols, called letters,  belonging to
a given set $\alpha$, called an alphabet. In these lectures we
discuss an approach to combinatorics of words based on their analogy
with curves on the plane. To begin with,  consider Figure \ref{fig1}
depicting a plane curve with   distinguished origin and
  orientation. The  three crossing points of the curve are
labeled by letters $a,b, c\in \alpha $. Now, starting at the origin
we move along the curve. The first crossing met is labeled $a$, the
second one   $b$, and so on. Finally we return to the origin and
stop. Writing down the labels of the crossings as we encounter them,
we obtain the word $abcabc$. This procedure, deriving a word from a
closed plane curve with labeled double transversal intersections
  was   introduced by Gauss \cite{ga}  in his
attempt to classify such curves. Clearly, every letter appears in
the resulting word  twice. Words in which every letter appears twice
are called Gauss words.

\begin{figure}[htbp]
\psfrag{a}[cc][cc]{$a$}
\psfrag{b}[cc][cc]{$b$}
\includegraphics[scale=0.6]{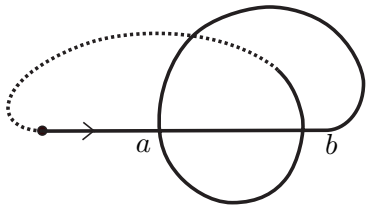}
\caption{The word $abab$ is not realizable by a closed curve on
$\R^2$} \label{fig2}
\end{figure}

It is easy to see that not all Gauss words can be realized by closed
plane curves. The word $abab$ for instance is   not realizable by
such a curve--- see Figure \ref{fig2}. Various conditions on  Gauss
words necessary and sufficient for their realizability by closed
plane curves were obtained by several authors, see   \cite{Ma},
\cite{LM}, \cite{ro}, \cite{dt}, \cite{cw}, \cite{ce}, \cite{cr}.

\begin{figure}[htbp]
\psfrag{a}[cc][cc]{$a$}
\psfrag{b}[cc][cc]{$b$}
\includegraphics[scale=0.6]{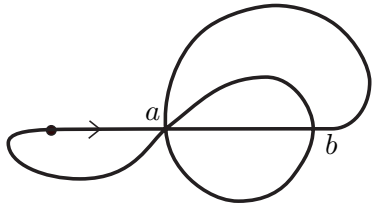}
\caption{A curve with a triple point yielding the word $ababa$}
\label{fig3}
\end{figure}

The aim of these lectures is  to study {\it arbitrary} words using
ideas taken from the topology of curves. To do this, we generalize
the above picture in the following three directions. First of all,
we allow curves with self-crossings of arbitrary multiplicity $\geq
2$. For example, the word associated with the curve in Figure
\ref{fig3} is $ababa $. In general, the number of occurrences of the
label of a crossing in the associated word is equal to the
multiplicity of this crossing.  To handle letters appearing only
once, we may distinguish a finite set of generic points on the curve
and label them as well; we shall not do that here.

\begin{figure}[htbp]
\psfrag{a}[cc][cc]{$a$}
\psfrag{b}[cc][cc]{$b$}
\psfrag{c}[cc][cc]{virtual crossing}
\includegraphics[scale=0.6]{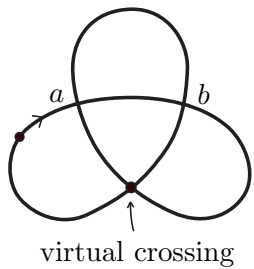}
\caption{A curve with an unlabeled crossing}
\label{fig4}
\end{figure}

The second direction in which we may generalize is to allow
\emph{unlabeled} or \emph{virtual} crossings. Such   crossings do
not contribute at all to the associated word. For example, the curve
in Figure \ref{fig4} gives rise to the word $abab $. The idea of
unlabeled crossings is inspired by the theory of  virtual knots
introduced by L. Kauffman \cite{ka}.

\begin{figure}[htbp]
\psfrag{a}[cc][cc]{$A_{1}$}
\psfrag{b}[cc][cc]{$A_{2}$}
\psfrag{c}[cc][cc]{$A_{3}$}
\psfrag{d}[cc][cc]{$a$}
\includegraphics[scale=0.65]{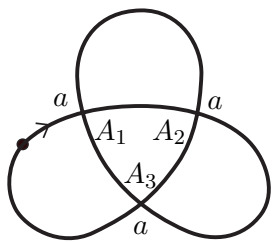}
\caption{A curve with   label $a$ on different crossings}
\label{fig5}
\end{figure}

Thirdly, we may allow the same label to be used on different
crossings. In Figure \ref{fig5} there are three crossings $A_{1} $,
$A_{2}$,   $A_{3} $,   all   labeled by the same letter $a \in
\alpha$. This   leads us to so-called \'etale words that are words
in an alphabet projecting to the given (fixed) alphabet $\alpha$.
The curve  on Figure \ref{fig5} gives rise to the \'etale word $A_1
A_2A_3 A_1A_2A_3$ where the letters $A_1, A_2, A_3$ project to $a\in
\alpha$. This \'etale word is
  a Gauss word in the alphabet $\{A_1, A_2, A_3\}$; we call such \'etale
words {\it nanowords}.

The general scheme  of the topology of words is as follows:
arbitrary words on the given alphabet $\alpha$  are approximated by
nanowords and the latter are studied by methods inspired by the
topology of curves. It is certainly interesting for topologists to
apply topological methods to
  study   such new objects.  Our approach also leads to new
  questions  concerning   combinatorics of words.
   The   accent in this theory shifts from words themselves to
transformations of words, inspired by topology. The situation is similar
to the one in
  knot theory where one focuses on isotopy classes of knots rather than on
  specific knot diagrams.

The present exposition follows my lectures  given in the Research
Institute for Mathematical Sciences (RIMS, Kyoto) in February 2006
and is based on my papers \cite{Tu0}--\cite{Tu3}. The lectures were
organized by Prof.\ Tomotada Ohtsuki. The  lecture notes taken by
Eri Hatakenaka, Daniel Moskovich and  Tadayuki Watanabe served as
the basis for this paper. I  would like to express my gratitude to
Prof.\ Ohtsuki for inviting me to RIMS and for organizing the
lectures and to Eri Hatakenaka, Daniel Moskovich and Tadayuki
Watanabe for   the
  preparation of the notes.

\section{Words and  nanowords}

In this section we give formal definitions  of words, \'etale words,  and
nanowords.  Fix a set $\alpha$ called   the
  {\it alphabet}.

\subsection{Words}
A {\it word of length} $n \geq 1$  on $\alpha$
is   a mapping
$$w\co\widehat {n} \to \alpha,$$
where $\widehat {n}$ is the set $\{1, 2,   \ldots, n\}.$
That is, a word of length $n$  on $\alpha$
is a sequence of $n$ elements of $\alpha $.
For example, the word $w=aba $,
with $a,$ $b\in \alpha,$
is nothing but the map
$$w\co\{1, 2, 3\}\to \alpha,$$
defined by
$w(1)=a,$ $w(2)=b$, and $w(3)=a.$
By convention, there is one empty word $\phi$ of length $0.$

The   {\it opposite word} to a word $w\co\widehat {n}\to \alpha,$ is
denoted by $w^{-} $ and defined by   $w^{-}(i)=w(n+1-i) $ for all
$i\in \widehat  n$. For instance, if $w=abca,$  then $w^{-}=acba.$

{\it Concatenation} of two words is   defined   by writing down
the first word and then the second one.
For example, the concatenation  of the   words $w=abc$ and $v=dbba$
on $\alpha $
is the word $wv=abcdbba$.

One more operation on words is a change of the alphabet. For a map
$f\co\alpha \to \alpha'$ from an alphabet $\alpha$ to an  alphabet
$\alpha' $ and a word $w\co\widehat {n} \to \alpha,$
  set $f_{\sharp}(w)=f\circ w$. This is a word  on $\alpha' $ of the same
length $n$.
For example, if $w=abca,$ then $f_{\sharp} (w)=f(a)f(b)f(c)f(a).$

\subsection{\'Etale words}
  An {\it $\alpha$--alphabet} is
  a set $\mathcal{A}$ endowed with a mapping to
$  \alpha,$ called the {\it projection}. The image in $\alpha$
of any  $A\in \mathcal{A}$ will be denoted $|A|.$

A {\it morphism} of $\alpha$--alphabets $\mathcal{A}_1$ and
$\mathcal{A}_2$ is a mapping $f\co\mathcal{A}_1\to
\mathcal{A}_2$ such that $|A|=|f(A)|$ for any  $A\in
\mathcal{A}_1.$ This means that the diagram
\[\xymatrix{
\mathcal{A}_1 \ar[rr]^{f} \ar[dr]_{\text{proj}} &
& \mathcal{A}_2 \ar[ld]^{\text{proj}}\\
& \alpha & \\
}\] commutes. An {\it isomorphism} of $\alpha$--alphabets is a
bijective morphism.

An {\it \'etale word} over $\alpha$ is   a pair (an
$\alpha$--alphabet $\mathcal{A}$, a word in the $\alpha$--alphabet
$\mathcal{A}$). In particular, every word on
$\alpha$ becomes an \'etale word over $\alpha$ by regarding $\alpha$ as an
$\alpha$--alphabet $\mathcal{A}=\alpha$ with   projection
to $\alpha$ being the identity map.
Another example: let $a\in  \alpha  $ and
$\mathcal{A}=\{A_1, A_2, A_3\}$   with      $|A_1|=|A_2|=|A_3|=a$.
The pair $(\mathcal{A}, A_1 A_2 A_3 A_1 A_2 A_3)$
is an \'etale word over $\alpha.$ It corresponds to
the picture  on Figure \ref{fig5}.

Two \'etale words $(\mathcal{A}_1, w_1)$ and $(\mathcal{A}_2, w_2)$
over   $\alpha$ are  {\it isomorphic} if there is an
isomorphism $f\co\mathcal{A}_1\to\mathcal{A}_2$ such that
$w_2=f_{\sharp}(w_1).$ The relation  of isomorphism will be denoted
$ \approx .$ For example, if   $\mathcal{B}=\{B_1, B_2, B_3\} $ is an
$\alpha$-alphabet with      $|B_1|=|B_2|=|B_3|=a\in \alpha$, then
$$(\mathcal{B}, B_1 B_2 B_3 B_1 B_2 B_3)\approx
(\mathcal{A}, A_1 A_2 A_3 A_1 A_2 A_3)$$ where the \'etale word on the right
is as in
the previous paragraph.

  For an \'etale word
$(\mathcal{A}, w)$, the {\it opposite \'etale word}
is defined by
$(\mathcal{A}, w)^{-}=(\mathcal{A}, w^-).$
The {\it product} of   \'etale words
$(\mathcal{A}_1, w_1)$ and $(\mathcal{A}_2, w_2)$
over   $\alpha$
is defined as follows. If $\mathcal{A}_1 \cap\mathcal{A}_2=\phi$,
then the pair
  $(\mathcal{A}_1  \cup\mathcal{A}_2, w_1 w_2) $
   is an \'etale word over $\alpha,$ and we call it
  the product of
  $(\mathcal{A}_1, w_1)$ and $(\mathcal{A}_2, w_2)$.
  If $\mathcal{A}_1 \cap\mathcal{A}_2\neq \phi$, then we
pick   an \'etale word
  $(\mathcal{A}_1', w_1')$ over $\alpha$ isomorphic to  $(\mathcal{A}_1, w_1)$
  and such that
  $\mathcal{A}_1' \cap\mathcal{A}_2=\phi$. We call then
  the  \'etale word
  $(\mathcal{A}_1' \cup\mathcal{A}_2, w_1'w_2),$
    the product of
  $(\mathcal{A}_1, w_1)$ and $(\mathcal{A}_2, w_2).$
The product of \'etale words is well defined up to
isomorphism.

Beware that   concatenation of words on
$\alpha$ usually differs from multiplication of the corresponding
\'etale words. For example, for the  words $w_1=abb$ and $w_2=aa$ on the
alphabet $\alpha=\{a, b\} $, the corresponding \'etale words are
$(\mathcal{A}_1=\{A,
B\}, ABB)$    with   $|A|=a, |B|=b,$ and   $(\mathcal{A}_2=\{A', B'\}, A'A')$
     with  $|A'|=a,  |B'|=b$. Their product is the
       \'etale word   $(\{A, B, A',B'\}, ABBA'A'),$  with   $|A|=|A'|=a, |B|=|B'|=b,$
       while the \'etale word corresponding to $w_1 w_2=abbaa$
is  $(\{A, B \}, ABBAA)$ with   $|A|=a, |B|=b$. These two \'etale words are not
isomorphic.

\subsection{Nanowords}
A word $w$ on a finite  alphabet   is   a {\it Gauss word} if every letter of this alphabet
appears in
$w$ exactly two times. For instance, the words
$   AABB, ABAB, ABBA, BAAB, BABA, BBAA,    $ are all Gauss words on the alphabet
$ \{A,B\}$.
The words $AA$, $ ABA$ are not
Gauss words on this alphabet.

An \'etale word $(\mathcal{A}, w)$ is   a {\it nanoword}
over $\alpha$ if $w$ is a Gauss word on $\mathcal{A}.$ Then
the alphabet $\mathcal{A}$ is
finite and the number of its elements is equal to the    half of the length of
$w$. For example, let $\mathcal{A}=\{A, B\}$
   with  $|A|=a\in
\alpha$ and $|B|=b\in \alpha$. The \'etale word $(\mathcal{A},
ABAB)$ is a nanoword. Another example: $\mathcal{A}=\{A, B,
C\}$   with $|A|=a\in \alpha, $ $|B|=b\in
\alpha$, and $|C|=c\in \alpha.$ Then the \'etale word $(\mathcal{A},
ABCBCA)$ is a nanoword.

Two nanowords $(\mathcal{A}_1, w_1)$ and $(\mathcal{A}_2, w_2)$ are
  {\it isomorphic} if they are isomorphic as \'etale words.

We now  make a few simple remarks about nanowords. If
$(\mathcal{A}, w)$ is a nanoword then its opposite $(\mathcal{A},
w)^-$ is also a nanoword. The concatenation of two nanowords is a
   nanoword. An empty \'etale word $ (\mathcal{A}=\phi,
w=\phi) $   is a nanoword.
The set of nanowords over $\alpha$ is infinite provided $\alpha\not=\phi.$

Note finally that a plane curve with     fixed origin and
labeled crossings gives rise to a
  nanoword if
and only if this curve is generic--- that is all its
self-intersections are double transverse crossings. Thus,  the
nanowords can be thought of as combinatorial analogues
  of
generic curves.

\subsection{Desingularization}

\begin{figure}[htbp]
\psfrag{a}[cc][cc]{$a$} \psfrag{b}[cc][cc]{$b$}
\psfrag{c}[br][bc]{$\ A_{1,2}$} \psfrag{d}[t][cc]{$A_{2,3}$}
\psfrag{e}[tl][cc]{$A_{1,3}$} \psfrag{f}[cc][cc]{$B$}
\includegraphics[scale=0.6]{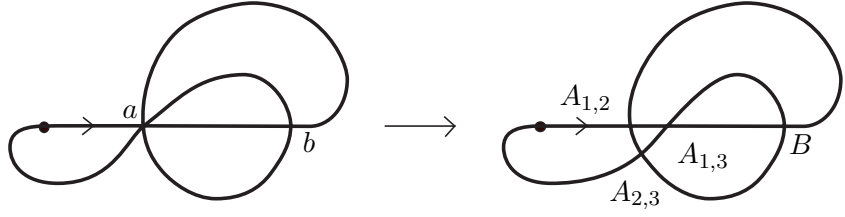}
\caption{Desingularization of a curve}
\label{fig6}
\end{figure}

Consider again the curve on Figure \ref{fig3}. By a small
deformation near the triple point, we obtain a generic curve whose
singularities are double points $A_{1,2}, A_{2,3}, A_{1,3}, B$ shown
on Figure \ref{fig6}.   We label the  points $A_{1,2}, A_{2,3},
A_{1,3}$ by   $a$ and the point $B$ by   $b$. Thus, each crossing of
the deformed curve has the same label as the corresponding crossing
of the original curve.
  The   left curve on Figure \ref{fig6} gives the word $w=ababa.$ The right curve
  gives the  Gauss word $w^d=A_{1,2} A_{1,3} B A_{1,2} A_{2,3} B A_{1,3} A_{2,3} $
   in the
  $\alpha$--alphabet $  \{A_{1,2}, A_{2,3}, A_3, B\}$   with
  $\abs{A_{1,2}}=\abs{A_{1,3}}=\abs{A_{2,3}}=a$ and $\abs{B}=b$.  We view
  the nanoword  $w^d$  over $\alpha$ as the desingularization of the word $w=ababa$.

The desingularization  procedure considered in this example can be generalized and
leads thus  to  desingularization of arbitrary \'etale words.
More precisely, for every \'etale word $(\mathcal{A}, w)$ over $\alpha,$ we define a
nanoword $(\mathcal{A}^d, w^d)$ over $\alpha$ called its
{\it desingularization}. For any   $A\in \mathcal{A},$ the {\it
multiplicity} $m_w (A)$ is   the number of times that
$A$ occurs in $w.$ For instance, $m_{ABB}(A)=1 $ and
$m_{ABB}(B)=2.$ We define the  set $\mathcal{A}^d$ by
$$\mathcal{A}^d=\{(A, i, j)\,|\,\, A\in \mathcal{A}, \,\,1\leq i<j\leq m_w (A)\}.$$
Denoting $(A, i, j)$ by $A_{i, j}$, we define the projection
$\mathcal{A}^d\to \alpha$ by $|A_{i, j}|=|A|\in \alpha.$ This makes ${\mathcal A}^d$
into an $\alpha$--alphabet. Here every letter of multiplicity $m \geq 1$
in   $w$ gives rise to
$\frac{m(m-1)}{2}$ letters in   ${\mathcal A}^d $.  The word
$w^d$ on the alphabet $\mathcal{A}^d$ is defined in   two steps.
\begin{enumerate}[Step 1.]
  \item Delete from $w$ all letters of multiplicity 1.
  \item For each $A\in \mathcal{A}$ with $m_w (A)\geq 2$
  and for each $i=1,~2,~\ldots,~m_w (A),$
  we replace the $i$-th entry of $A$ in $w$ by
  $$A_{1, i}~A_{2, i}~\cdots~A_{i-1, i}~A_{i, i+1}~A_{i, i+2}~\cdots~A_{i, m_w (A)}.$$
\end{enumerate}
The resulting  word $w^d$ on  $\mathcal{A}^d $   is a
Gauss word and the
  pair $(\mathcal{A}^d, w^d)$ is   a nanoword over $\alpha$. For
example, for the \'etale word  $(\mathcal{A}=\{A, B\}, w=ABABA)$   with
$|A|=a, |B|=b $, this procedure gives  the nanoword
$$(\mathcal{A}^d, w^d)=(\{A_{1,2}, A_{2,3}, A_{1,3}, B_{1,2}\},
A_{1,2}A_{1,3}B_{1,2}A_{1,2}A_{2,3}B_{1,2}A_{1,3}A_{2,3})$$
with $|A_{1,2}|= |A_{2,3}|= |A_{1,3}|=a$ and $|B_{1,2}|=b$.

  We can say   that the
desingularization of curves  is based on  viewing the crossings under
a   strong microscope which allows us to see the \lq \lq internal structure" of
each
crossing. This analogy suggested the term    nanoword.

\section{Knot theory and homotopy of nanowords}

\subsection{Knot theory}
\begin{figure}[htbp]
\psfrag{a}[cc][cc]{$1$}
\psfrag{b}[cc][cc]{$2$}
\psfrag{c}[cc][cc]{$3$}
\includegraphics[scale=0.55]{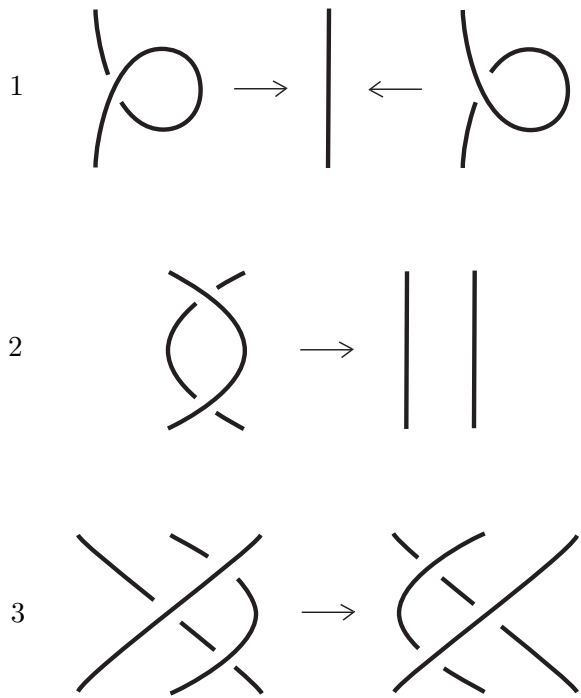}
\caption{Reidemeister moves on knot diagrams}
\label{fig8}
\end{figure}

We shall study nanowords using the analogy with curves and knots.
Recall the   classical Reidemeister moves on knot diagrams in
$\R^2$,   shown in Figure \ref{fig8}. (The inverse moves are also
called the Reidemeister moves.) These moves  and the isotopy of knot
diagrams in $\R^2$  generate an equivalence relation on   knot
diagrams which we call   the {\it R--equivalence}. We have
$$\{\text{knots in $\R^3$}\}/\text{ isotopy}
=\{\text{knot diagrams}\}/\text{ R--equivalence}.$$ This fundamental
equality, due to K.\ Reidemeister,   reduces the study of isotopy
classes of knots to a study of R--equivalence classes of knot
diagrams.

Similar moves can be considered on pointed  curves. For the sake of
the following discussion, we switch to curves--- that is we make no
distinction  between under-crossings and over-crossings. We always
assume that the moves act away from the origins of the curves.

\begin{figure}[htbp]
\psfrag{a}[cc][cc]{$A$}
\includegraphics[scale=0.6]{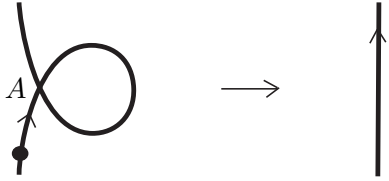}
\caption{First Reidemeister move  on a labeled curve}
\label{fig9}
\end{figure}

Let us look at the effect of the Reidemeister  moves on words
associated with curves. The first Reidemeister move, shown in Figure
\ref{fig9}, acts   as  $xAAy\mapsto xy,$ where $x$ and $y$ are words
not including the letter $A$.

\begin{figure}[htbp]
\psfrag{a}[cc][cc]{$A$}
\psfrag{b}[cc][cc]{$B$}
\includegraphics[scale=0.6]{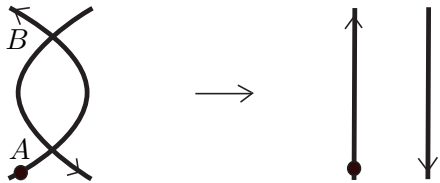}
\caption{Second Reidemeister move  on a labeled curve}
\label{fig10}
\end{figure}

Consider the second Reidemeister move  with labels, orientations, and the position
of the origin as   in Figure \ref{fig10}.    The move  acts on the associated
word  as
$xAByBAz\mapsto xyz $ where $x, y,z$ are words not including the letters
$A, B$.  For another choice  of orientations, the move may act
as
$xAByABz\mapsto xyz $.   The first version   $xAByBAz\mapsto xyz $ is sufficient for our
  aims as will be clear from the results below.

\begin{figure}[htbp]
\psfrag{a}[cc][cc]{$A$}
\psfrag{b}[cc][cc]{$B$}
\psfrag{c}[cc][cc]{$C$}
\includegraphics[scale=0.6]{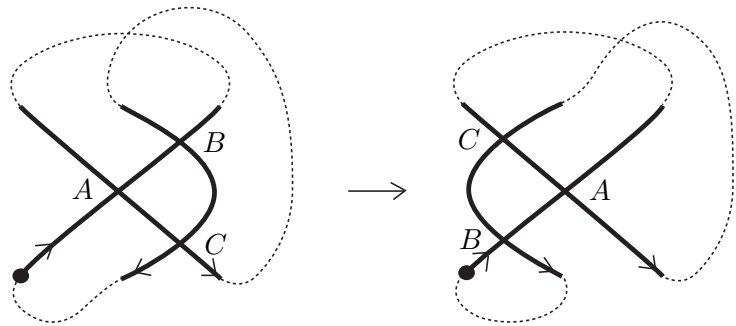}
\caption{Third Reidemeister move on a labeled curve}
\label{fig11}
\end{figure}

Consider the third
Reidemeister move  with labels, orientations, the position
of the origin, and the order of branches  as   in Figure \ref{fig11}.
This move  acts on the associated
word  as   $xAByACzBCt\mapsto xBAyCAzCBt $ where $x, y,z, t$
  are words not including the letters
$A, B, C$. For  other choices  of orientations, order of branches etc.,
the move may act
differently but   the version  shown on Figure \ref{fig11}
is sufficient for our
  aims.

\subsection{Homotopy of nanowords}  Let $\alpha$ be an alphabet (a fixed set).
We fix {\it homotopy data} consisting of an involution $\tau\co\alpha\to
\alpha$ and an arbitrary  set ${\mathcal S}\subset
\alpha^3=\alpha\times\alpha\times\alpha.$  The geometric meaning of  $\tau$ and
${\mathcal S}$ will be discussed in the next
section. In the context of knot diagrams,   $\tau$   switches between positive and negative
crossings.  The role of ${\mathcal S}$ is to determine the    transformations of
labels accompanying the third
Reidemeister move.

Three {\it homotopy moves}
on nanowords over $\alpha$ are defined as follows.
\begin{enumerate}[(1)]
  \item $(\mathcal{A}, xAAy)\mapsto (\mathcal{A}-\{A\}, xy) $ where $x, y$ are words on the alphabet
  $\mathcal A- \{A\}$.
  Note that if $(\mathcal{A}, xAAy)$ is a nanoword, then so is $(\mathcal{A}-\{A\}, xy) $.
  The inverse move   adds a new letter $A$ to the
  $\alpha$--alphabet and   inserts $AA$ into the word.
  \item $(\mathcal{A}, xAByBAz)\mapsto (\mathcal{A}-\{A, B\}, xyz)$
provided  $|A|=\tau(|B|)$ and
  $x,$ $y$, $z$ are words on the alphabet $\mathcal{A}-\{A, B\}.$
  \item $(\mathcal{A}, xAByACzBCz)\mapsto (\mathcal{A}, xBAyCAzCBt),$
provided $(|A|, |B|, |C|)\in {\mathcal S}  $ and
  $x,$ $y$, $z$, $t$ are words on the alphabet $\mathcal{A}-\{A, B, C\}.$
\end{enumerate}

Two nanowords over $\alpha$ are   {\it ${\mathcal S}$--homotopic} if
they can be related by a finite sequence of homotopy moves, inverse
moves,  and isomorphisms. We denote this equivalence relation by
$\simeq_{\mathcal S}  $ and call it {\it ${\mathcal S}$--homotopy.}
This definition readily extends to \'etale words:
  \'etale words $w_1$ and $w_2$ are   {\it ${\mathcal S}$--homotopic}
  if $w_1^d\simeq_{\mathcal S} w_2^d$.
  In particular, the notion of ${\mathcal S}$--homotopy applies
  to   words on $\alpha$.

  We will use the following notation:
$$
\mathcal{N}(\alpha,{\mathcal S})=\set{\text{set of nanowords over
$\alpha$}}\big/{\mathcal S}\text{--homotopy}.
$$
Clearly,
$\mathcal{N}(\alpha,{\mathcal S})$ is a monoid, with the empty nanoword
as its unit
element and concatenation as its product. This monoid   depends on $\tau$ which is however omitted
in the notation $
\mathcal{N}(\alpha,{\mathcal S})$  to make it shorter.

The following two lemmas
show that our three moves generate a wider set of similar moves.
In the context of Figures \ref{fig10} and \ref{fig11}, the  new
moves
correspond to other
  choices of
orientations, branch connections, etc.

\begin{Lem}  Let $A,B, C$ be three distinct
letters in an $\alpha$--alphabet $\mathcal A$ and let $x, y,z,t$
be words in the alphabet $\mathcal A- \{A, B, C\}$ such that $xyzt$ is a Gauss word in this
alphabet. Then
we have the following   ${\mathcal S}$--equivalences:
\begin{equation}
(\mathcal{A}, xAByCAzBCt)\simeq_{\mathcal S}(\mathcal{A}, xBAyACzCBt)
\tag{1}\end{equation}
if $(|A|, \tau(|B|), |C|)\in {\mathcal S},$
\begin{equation}
(\mathcal{A}, xAByCAzCBt)\simeq_{\mathcal S}(\mathcal{A}, xBAyACzBCt)
\tag{2}\end{equation}
if $(\tau(|A|), \tau(|B|), |C|)\in {\mathcal S},$ and
\begin{equation}
(\mathcal{A}, xAByACzCBt)\simeq_{\mathcal S}(\mathcal{A}, xBAyCAzBCt)
\tag{3}\end{equation}
if $(|A|, \tau(|B|), \tau(|C|))\in {\mathcal S}.$
\label{l1}
\end{Lem}

\begin{Lem}
Suppose that ${\mathcal S}\cap(\alpha\times b\times b)\not= \phi$
for all $b\in \alpha.$ Let $(\mathcal A, xAByABz)$ be a nanoword over $\alpha$ where $A, B\in
\mathcal A$ with $|A|=\tau(|B|)$ and $x,y,z$ are words on the alphabet $\mathcal A- \{A,B\}$.
Then
$$(\mathcal{A}, xAByABz)\simeq_{\mathcal S}
(\mathcal{A}-\{A, B\}, xyz).$$
\label{L2}\end{Lem}
\begin{proof}
Set $b=|B|\in \alpha.$ By assumption, there is   $e\in
\alpha$ such that $(e, b, b)\in {\mathcal S}.$ Pick a letter $E$ not belonging to $\mathcal A $
and set $|E|=\tau(e)\in \alpha.$ Then
\begin{align*}
(\mathcal{A}, xAByABz)\overset{(\text{Move 1})^{-1}}{\simeq_{\mathcal S}}
&(\mathcal{A}\cup \{E\}, xAEEByABz)\\
\simeq_{\mathcal S}\hspace{15pt}&(\mathcal{A}\cup\{E\}, xEABEyBAz)\\
\overset{(\text{Move 2})}{\simeq_{\mathcal S}}
\hspace{5pt}&(\mathcal{A}\cup\{E\}-\{A, B\}, xEEyz)\\
\overset{(\text{Move 1})}{\simeq_{\mathcal S}}\hspace{5pt}
&(\mathcal{A}-\{A, B\}, xyz)
\end{align*}
In the second line we use the ${\mathcal S}$--homotopy   of Lemma \ref{l1}.(2), where
$A$ is replaced by $E,$ $B$ by $A,$ and $C$ by $B$. This homotopy applies
  since
$(\tau(|E|), \tau(|A|), |B|)$ $=(e, b, b)\in {\mathcal S}.$
\end{proof}

\subsection{Typical questions}
In analogy with knot theory, the main  objective of the homotopy theory of words
is to  classify \'etale words and nanowords up to  ${\mathcal S}$--homotopy.
Putting it differently, the goal is to compute the monoid $\mathcal{N}(\alpha,{\mathcal S})$ at least for
some choices of $\alpha$, $\tau$, and ${\mathcal S}$.
We are
very far from reaching this goal. Available results are outlined in
  the rest of the paper.

Taking knot theory as a model, we list here
  several  typical questions concerning the homotopy of words.

\begin{enumerate}[(Q1)]
  \item Is a given nanoword ${\mathcal S}$--contractible, {\it i.e.},  ${\mathcal S}$--homotopic to the
   empty
  nanoword ?
  This question corresponds to the question of whether or not a given knot diagram
presents an unknot.
  \item Is a given nanoword $w$ homotopically symmetric, that is
   ${\mathcal S}$--homotopic to $w^{-}$ ?
  Note that opposite words corresponds to knots with opposite orientations.
  \item    Define the {\it length norm} of an \'etale  word $w$ by
  $$||w||_{{\mathcal S}}=\frac{1}{2}~
  (\text{minimal length of a nanoword ${\mathcal S}$-homotopic to $w$}).$$
  Note that $||w||_{\mathcal S}=0$ if and only if $w$ is contractible.
  For any \'etale words $w_1, w_2$,
  $$||w_1 w_2||_{\mathcal S}\leq ||w_1||_{\mathcal S}+||w_2||_{\mathcal S}.$$
  Compute the length norm.
\end{enumerate}






\section{Curves and knots   as nanowords}

In this section we   clarify the relations  between curves, knots,
and nanowords.

\subsection{Curves    as nanowords}
\begin{figure}[htbp]
\psfrag{base point}[cc][cc]{basepoint}
\psfrag{regular\r}[cc][cc]{regular} \psfrag{neighborhood}
[cc][bc]{neighborhood}
\includegraphics[scale=0.6]{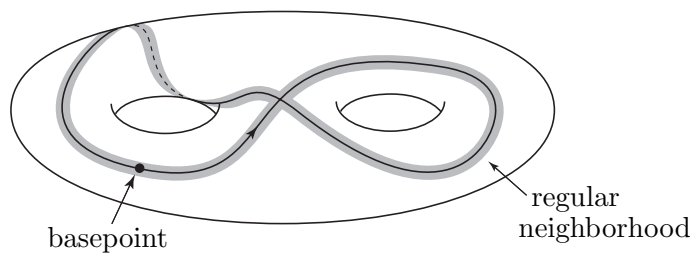}
\caption{A curve on a surface of genus two} \label{2-1}
\end{figure}

In the sequel, the word \lq\lq curve" means the image of a generic
immersion of an oriented circle into an oriented surface. Here
``generic'' means that the curve has only a finite set of
self-intersections, which are all double and transversal. The curve
may be immersed into any oriented surface, of any genus, compact or
not, with boundary or not. Note that all self-intersections of a
curve look locally like
$\begin{minipage}{10pt}\includegraphics[width=10pt]{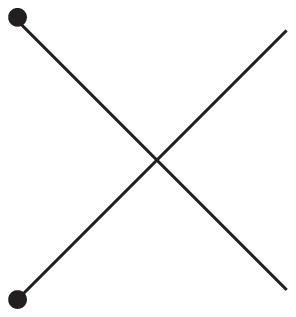}\end{minipage}$\hspace{5pt}.
Triple points and   self-tangencies are not allowed. Every curve has
a \emph{regular neighborhood}. This  is a narrow band around the
curve inside the surface, see Figure \ref{2-1}. Note that the
orientation of the ambient surface induces an orientation of the
regular neighborhood. \par

A curve is   \emph{pointed} if it is provided with a
base-point which is not a self-intersection.
An example of a pointed curve is drawn on Figure \ref{2-1}:\par

Two pointed curves are   \emph{stably homeomorphic} if their  regular neighborhoods
look exactly the same including the position of the curves in these neighborhoods. Here is a
more
precise definition.

\begin{defn}
Two pointed curves are     stably homeomorphic  if
there is an orientation-preserving homeomorphism of their regular neighborhoods mapping the
first curve onto the second one and preserving the origin and the
orientations of the curves.
\end{defn}

The stable homeomorphism class of a curve
is determined solely by the germ of the ambient surface near the
curve;  what happens outside a regular neighborhood  does not matter.
In particular, adding handles and punctures to the surface away from a
  neighborhood of the  curve does not change its stable homeomorphism
class.

Recall the definition of the \emph{stable equivalence} of curves from    \cite{kk}, \cite{cks}.

\begin{defn}
Two pointed curves are    stably equivalent  if they
can be related by a finite sequence  of the following transformations:
\begin{enumerate}
\item Stable homeomorphism.
\item Homotopy of the curve in its ambient surface away from the
origin.
\end{enumerate}
\end{defn}

The   homotopy in (2) may push a branch of the curve across another
branch or across a double point but not across the origin of the curve.

Pointed curves related by any sequence of   moves (1), (2) are
stably equivalent. Thus,  we may start with a curve, transform it by
a stable homeomorphism, deform the resulting curve, add handles,
deform again, puncture the surface,  etc. All these transformations
preserve the stable equivalence class of the curve.   As an
exercise, the reader may show that any two pointed curves on the
2--sphere are stably equivalent. The same is true for curves on the
2--torus. Pointed curves on surfaces of higher genus are not
necessarily stably equivalent. The classification of   stable
equivalence classes of pointed curves is an interesting topological
problem.

We note here three geometric invariants
  of
  pointed  curves preserved under the stable equivalence:
  the  minimal crossing number, the genus, and the virtual number.
  The minimal crossing number  $\Vert c\Vert$ of a  pointed  curve $c$ is the minimal
  number of crossings of a pointed curve stably equivalent to $c$.
  The genus   $g(c)$   is the
minimal integer $g\geq 0$ such that $c$ is stably equivalent to a pointed
curve on a  closed
surface of genus $g$. To define the virtual number of $c$, note that any
  pointed curve on $\R^2$
with a distinguished set of \lq\lq virtual" crossings represent a stable equivalence
class of pointed curves. One simply does not look at the   virtual crossings or,
equivalently,
trades a branch of $c$ near each virtual
crossing for a branch going along a small
1-handle attached to $\R^2$ and avoiding the rest of the curve.
The virtual number $v(c)$ is the
minimal integer $v\geq 0$ such that there is a pointed curve on $\R^2$
with $v$ virtual crossings representing the stable equivalence class of $c$.
It is clear that $v(c)\geq g(c)$. \par

Denote by $\mathcal{C}$  the set of stable
equivalence classes of pointed curves. The elements of $\mathcal{C}$ are
called   long flat knots \cite{ka} or open virtual strings \cite{Tu0}.\par

We now relate the theory of curves with the theory of nanowords.
Consider the following homotopy data:
$$
\alpha_0 =\{a,b\},\,\, \tau_0(a)=b, \tau_0(b)=a, \,\,\,
{\mathcal S}_0=\set{(a,a,a),(b,b,b)}\subset\alpha_0^{3}.
$$

\begin{Thm}
There is a canonical bijection
$$
\mathcal{C}\overset{\approx}{\To} \mathcal{N}(\alpha_{0},{\mathcal S}_{0}).
$$
\end{Thm}

\begin{figure}[htbp]
\centering
\psfrag{A}[cc][cc]{$A_{i}$}\psfrag{i}{}\psfrag{1-st}[cc][bc]{$1$st}
\includegraphics[scale=0.6]{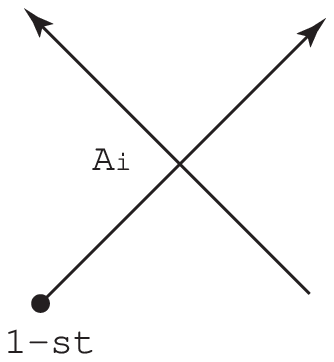}\hspace{1in}
\includegraphics[scale=0.6]{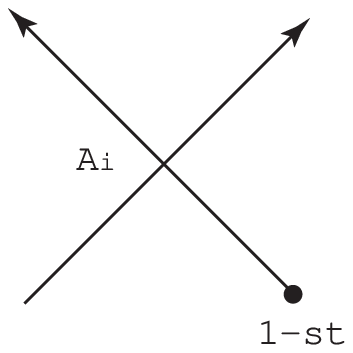}
\caption{On the left $\abs{A_{i}}=a$ and   on the right $\abs{A_{i}}=b$} \label{2-6}
\end{figure}

This theorem shows that the theory of nanowords includes the theory of pointed curves as a
special case.  We outline   a construction of the bijection $
\mathcal{C}\to \mathcal{N}(\alpha_{0},{\mathcal S}_{0})$.
Consider  a pointed curve on a surface. Label its crossings in an arbitrary way  by
different letters $A_{1},A_{2},\ldots,A_{n}$
where $n$ is the number of crossings. The Gauss word of the curve is obtained by
moving along the curve starting  at the
origin and writing down  the letters as we encounter them, finishing when
we get back to the origin. The resulting word, $w$,  on the alphabet
$$
\mathcal{A}=\set{A_{1},A_{2},\ldots,A_{n}}
$$
contains every letter
  $A_{1},A_{2},\ldots,A_{n}$ twice. We provide $\mathcal A$ with the projection to $\alpha_0$ as follows.
Consider the crossing of the curve  labeled $A_{i}$. If when
moving as above along the curve, we first traverse  this crossing
  from the bottom-left to the top-right, then
$\abs{A_{i}}=a$, otherwise   $\abs{A_{i}}=b$; see Figure \ref{2-6}, where
the orientation of the ambient surface is counterclockwise. The dot on the left (resp.\ right) picture is the
bottom-left (resp.\ bottom-right) entry of the crossing.
In this way  the set $\mathcal A$ becomes an $\alpha_0$--alphabet.
We assign to our curve the   class of this nanoword in
$\mathcal{N}(\alpha_{0},{\mathcal S}_{0})$.  We must prove that stably equivalent curves give rise to
${\mathcal S}_0$--homotopic nanowords. A different choice of the labeling of the crossings gives an
isomorphic nanoword. If   the curve is changed by a stable homeomorphism, then the
associated  nanoword does not
change, since it is defined entirely by the behavior of the curve  in its regular
neighborhood. A homotopy of the curve may be split into a composition of local
Reidemeister moves and the inverse moves. Then one verifies that under these moves the
associated nanoword  changes via the ${\mathcal S}_0$--homotopy moves
and the transformations   in Lemmas \ref{l1} and  \ref{L2}. The
resulting mapping  $
\mathcal{C}\to \mathcal{N}(\alpha_{0},{\mathcal S}_{0})$ is bijective, see \cite{Tu2}.  \par

Under the identification $
\mathcal{C}= \mathcal{N}(\alpha_{0},{\mathcal S}_{0})
$ the minimal crossing number of curves corresponds to the length
norm   on $\mathcal{N}(\alpha_{0},{\mathcal S}_{0})$.
The genus and the virtual  number yield interesting geometric
invariants of nanowords over $\alpha_0$.

One can suppress all references to the origin in the definitions above. This gives a
  relation of stable homotopy for non-pointed curves.  To obtain
  a corresponding notion for the
  nanowords, one has to introduce an additional move on nanowords, the so-called circular shift.
Briefly speaking, the shift moves the last letter of the word to the first position. For details,
see \cite{Tu2}.

\subsection{Knots  as nanowords}  
\begin{figure}[htbp]
\includegraphics[scale=0.6]{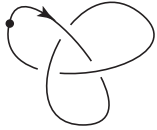}
\caption{A pointed knot diagram} \label{2-11}
\end{figure}

The constructions of the previous section can be
upgraded to the setting of knot diagrams.  By a
(pointed) knot diagram we mean a (pointed) curve on an oriented surface
with additional data at each crossing: one of the branches lies
``over'' and the other one lies ``under''.
A pointed knot diagram on $\R^2$ is shown on Figure \ref{2-11}.\par

Two pointed knot diagrams are  said to be  \emph{stably homeomorphic}
if there is an orientation-preserving
homeomorphism of the  regular neighborhoods  of the underlying curves, sending
the first    diagram onto the second one and preserving the
origin, the orientation, and the over-/under-crossing data.

Recall the   \emph{stable equivalence} of knot diagrams from    \cite{kk}, \cite{cks}.

\begin{defn}
Two pointed knot diagrams are    stably equivalent  if
they can be related by a finite sequence of the following transformations:
\begin{enumerate}
\item Stable homeomorphism.
\item The   Reidemeister moves on a
knot diagram in its ambient surface away from the origin.
\end{enumerate}
\end{defn}

The
moves  in (2) may push a branch of the diagram above or below a double
point or another branch but not across the origin.
Thus,   the origin may not lie inside   the neighborhoods
where the Reidemeister moves are   performed.\par

Let $\mathcal{K}$ denote the set of stable equivalence classes of
pointed knot diagrams. Elements of $\mathcal{K}$ are called
\emph{long virtual knots}, see \cite{ka}, \cite{gpv}.\par

The set $\mathcal{K}$ includes the set of isotopy classes of classical knots.
Classical knots are oriented knots in $S^{3}$. There is a map
$$
\set{\text{Classical knots}}/\text{isotopy} \hookrightarrow \mathcal{K}
$$
obtained by picking an arbitrary diagram of the given classical knot, picking an arbitrary
origin on the diagram (not a crossing), and taking the stable equivalence class of the
resulting pointed knot diagram.
This gives a well-defined mapping from the set of isotopy classes of
classical knots
into $  \mathcal{K}$. This mapping is known
to be injective, see \cite{ka}, \cite{gpv}.

The set $\mathcal{K}$ can be interpreted in terms  of nanowords as follows,
Set
\begin{itemize}
\item $\alpha_{*}=(a_{+},a_{-},b_{+},b_{-})$,
\item $\tau_{*}:\alpha_{*}\to \alpha_{*}$ the involution defined by
$\tau_{*}(a_{+})=b_{-}$ and $\tau_{*}(a_{-})=b_{+}$,
\item ${\mathcal S}_{\star}=\{(a_{\pm},a_{\pm},a_{\pm}),(a_{\pm},a_{\pm},a_{\mp}),
(a_{\mp},a_{\pm},a_{\pm}),$\\\hspace{20mm}
$(b_{\pm},b_{\pm},b_{\pm}),(b_{\pm},b_{\pm},b_{\mp}),(b_{\mp},b_{\pm},b_{\pm})\}
\subset\alpha_{*}^{3}$.
\end{itemize}

\begin{Thm}
There is a canonical bijection
$$
\mathcal{K}\overset{\approx}{\To}\mathcal{N}(\alpha_{*},{\mathcal S}_{*})
$$
such that  the following   diagram commutes:
\[
\begin{CD}
\mathcal{K}
@>\approx>> \mathcal{N}(\alpha_{*},{\mathcal S}_{*})\\
  @VVV @VVV\\
  \mathcal{C} @>\approx>> \mathcal{N}(\alpha_0,{\mathcal S}_0).
\end{CD}
\]
\noindent Here the map   $\mathcal{K}\to \mathcal{C}$ is given
by forgetting the over-/under-crossing data, and the map
$\mathcal{N}(\alpha_{*},{\mathcal S}_{*})\to \mathcal{N}(\alpha_0,{\mathcal S}_0)$ is given
by $a_{\pm}\mapsto a$, $ b_{\pm}\mapsto b$.
\end{Thm}

\begin{figure}[htbp]
\centering
\psfrag{A}[rr][ll]{$A_{i}$}\psfrag{i}{}\psfrag{1-st}[cc][bc]{$1$st}
\includegraphics[scale=0.5]{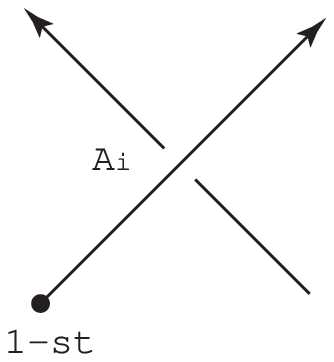}\hspace{0.6in}
\includegraphics[scale=0.5]{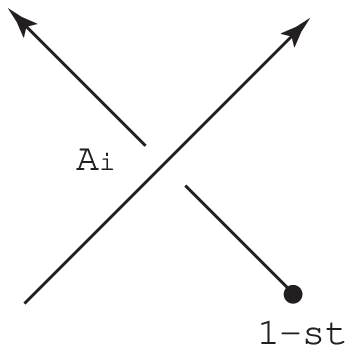}
\hspace{0.6in}
\includegraphics[scale=0.5]{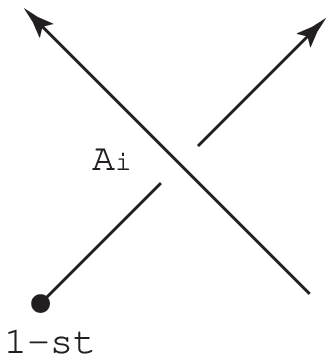}
\hspace{0.6in}
\includegraphics[scale=0.5]{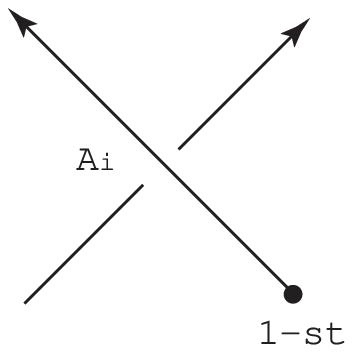}
\captionstyle{center} \onelinecaptionsfalse
\caption{From left to right:\protect\\
$\abs{A_{i}}=a_{+}$, $\abs{A_{i}}=b_{+}$, $\abs{A_{i}}=a_{-}$, and
$\abs{A_{i}}=b_{-}$} \label{2-77}
\end{figure}

This theorem shows that the theory of nanowords includes the theory of long virtual knots   as a
special case.
The definition of the bijection   $\mathcal{K}\to  \mathcal{N}(\alpha_{*},{\mathcal S}_{*})$ goes
similarly to the one for curves. The difference is that now  we project to $\alpha_\ast$ rather
than to $\alpha_0$. The rule is shown on Figure \ref{2-77}.

To extend these ideas to    links, one has to involve phrases, {\it i.e.},
sequences of words, see \cite{Tu2}, \cite{Tu3}.

\subsection{An extension to general $\alpha$}\label{sharrr} Let $\alpha$ be an
arbitrary alphabet with homotopy data $\tau, {\mathcal S}$.
Nanowords over   $\alpha$ can  be geometrically interpreted as
follows. Pick a  mapping $f\co\alpha\to \alpha_0$    such that
$f\tau=\tau_0 f$ and $f(a)=f(b)=f(c)$ for any triple $(a,b,c)\in
{\mathcal S}$.  Every nanoword $(\mathcal A, w\co\widehat n\to
\alpha)$ over $\alpha$ determines a nanoword $f_{\sharp}(w)=
(\mathcal A, f w\co\widehat n\to \alpha_0)$ over $\alpha_0$. The
latter can be represented by a pointed curve on a surface. Thus, $w$
gives rise to a family of  pointed curves $\{f_{\sharp}(w)\}_f$ {\it
underlying $w$} and numerated by   $f$ as above. Geometric
invariants of these curves   provide   geometric information about
$w$.  Stable equivalence classes of these curves depend  only on the
${\mathcal S}$--homotopy class of $w$.

This geometric representation of $w$ seems to be especially efficient
in the case where ${\mathcal S}$ is the
diagonal of $\alpha^3$ so that the conditions on $f$  reduce
  to the equivariance relation
$f\tau=\tau_0 f$. One interesting question: when does a given family
of stable equivalence classes of pointed curves numerated by
equivariant maps $\alpha\to \alpha_0$  arise from a
  nanoword over $\alpha$ ?

A similar geometric interpretation of
nanowords
in terms of knot
diagrams can be obtained by replacing $\alpha_0$ with
$\alpha_\ast$.





\section{Invariant $\gamma$ and self-linking   of nanowords}

In this section and in the sequel,
the symbol $\alpha$ denotes a (fixed) alphabet with involution
$\tau\co \alpha\to\alpha$.

\subsection{The set ${\mathcal S}$}

  For
the rest of this paper, we
set
$$
{\mathcal S}=\text{diagonal}=\set{(a,a,a)}_{a\in\alpha}
$$
The third homotopy move
is  $xAByACzBCt\mapsto xyz$ provided $\abs{A}=\abs{B}=\abs{C}$. In the sequel,
we   leave ${\mathcal S}$ out of notation. By \emph{homotopy} of nanowords,
   we   mean
${\mathcal S}$--homotopy when ${\mathcal S}$ is the diagonal as above. The homotopy
relation is denoted $\simeq$. The case of knots is   excluded
by this choice of ${\mathcal S}$, but the case of curves is covered. \par

We give now an example of homotopic words.  Pick $a\in\alpha$ such that $\tau(a)\neq a$.
Set $b=\tau(a)$. Consider the words:
$$w_{1} =aabab,\,\,
w_{2} =babaa,\,\,
w_{3} =baaab$$
and the nanoword
$$w_{4} = (\{A,A'\},AA'AA')\text{ where $\abs{A}=\abs{A'}=a$.  } $$

\begin{claim}
$w_{1}\simeq w_{2}\simeq w_{3}\simeq w_{4}$
\end{claim}

\begin{proof}
We prove only that $w_{1}\simeq w_{4}$. The proofs that $w_{2}\simeq
w_{4}$ and $w_{3}\simeq w_{4}$ are similar.\par
The desingularization of $w_1$ gives
$$
w_{1}^{d}=
A_{1,2}A_{1,3}A_{1,2}A_{2,3}BA_{1,3}A_{2,3}B
$$
where $\abs{A_{1,2}}=\abs{A_{1,3}}=\abs{A_{2,3}}=a$ and $\abs{B}=b$.
By Lemma \ref{L2}, we can strike out the two occurrences of
$A_{2,3}B$. This gives the nanoword $ A_{1,2}A_{1,3}A_{1,2}A_{1,3}$
isomorphic to $w_{4}$.
\end{proof}

This example   shows that the relation of homotopy  is quite non-trivial.

\subsection{A group--theoretic homotopy invariant}

  We   construct here a homotopy invariant of nanowords, $\gamma$.
  First, define a group $\Pi$ by generators and relations:
$$
\Pi=\langle\set{z_{a}}_{a\in\alpha}\mid z_{a}z_{\tau(a)}=1\text{ for
all $a\in\alpha$}\rangle.
$$
Note that if $a\neq\tau(a)$, then we have a free
generator $z_{a}=(z_{\tau (a)})^{-1}$ of $\Pi$, and if $a=\tau(a)$,
then  we have a generator $z_{a}$ of
order $2$.

For a nanoword $(\mathcal{A}, w\co\widehat {n}\to\mathcal{A})$ of length $n$, we
define $n$ elements $\gamma_{1}$, $\gamma_{2}$, $\ldots$, $\gamma_{n}$ of $\Pi$ by:
$$
\gamma_{i}=\left\{
             \begin{array}{ll}
               z_{\abs{w(i)}}, & \hbox{if $w(i)\neq w(j)$ for all $j<i$;} \\
               z^{-1}_{\abs{w(i)}}, & \hbox{otherwise.}
             \end{array}
           \right.
$$
Here $w(i)\in \mathcal A$ is the $i$-th letter of $w$ and
$\abs{w(i)}\in\alpha$ is its projection to $\alpha$.
The sequence $\gamma_{1},\gamma_{2},\ldots,\gamma_{n}$
may be also described as follows.
Since $w$ is a nanoword, each letter of $\mathcal A$
appears twice in the sequence $w(1), w(2), \ldots , w(n)$.
The first time it appears we write at this place the corresponding generator
of $\Pi$. At
its second appearance, we write down the inverse of that generator. This procedure
gives the sequence $\gamma_{1},\gamma_{2},\ldots,\gamma_{n}$.
Set
$$
\gamma(w)=\gamma_{1}\gamma_{2}\cdots\gamma_{n}\in\Pi.
$$
Since each generator appears in this product twice with opposite powers, the
abelianization of $\gamma(w)$ is zero. Thus, $\gamma(w)$ lies is in
the commutator subgroup $[\Pi,\Pi]\subset\Pi$.

For example, consider the nanoword $w=ABAB$ with
$\abs{A}=a\in\alpha$ and $\abs{B}=b\in\alpha$. Then
$\gamma(w)=z_{a}z_{b}z^{-1}_{a}z^{-1}_{b}\in [\Pi,\Pi]$.\par

\begin{Thm}
$\gamma(w)$ is a homotopy invariant of $w$.
\end{Thm}

\begin{proof}[Proof (outline)]  Under the first homotopy move,
$
\gamma(  xAAy )=\gamma(  xy )
$
because the first appearance of $A$ contributes $z_{a}$ and
the second appearance of $A$ contributes $z_{a}^{-1}$ where $a=\abs{A}$.
So we have the invariance under the
first homotopy move. The other two moves are treated similarly.
\end{proof}

The mapping
$$
\gamma\co \mathcal{N}(\alpha)=\mathcal{N}(\alpha,{\mathcal S})\To[\Pi,\Pi]
$$
is a monoid homomorphism (${\mathcal S}$ is   the diagonal).
It is easy to check that it is surjective.

The group $[\Pi,\Pi]$ can be shown to be free  for any $\tau$. If $\tau$ has at
least two orbits, then this group is   non-trivial
and by the results above,
  $\mathcal{N}(\alpha)$ is
an infinite monoid. If $\tau$ has at least three orbits, then $[\Pi,\Pi]$
has rank $\geq 2$,  and by the results above,
$\mathcal{N}(\alpha)$ is a non-abelian monoid.\par

Let us consider two examples
   where   $\gamma$ does not work.
The interesting case in topology is the case of curves, where
$\alpha=\set{a,b}$ with $\tau(a)=b$. In
this case:
$$
\Pi=\langle z_{a},z_{b}\mid z_{a}z_{b}=1\rangle=\Z
$$
and  $
[\Pi,\Pi]=0
$. So, for topology the invariant $\gamma$ is of no interest.
Another example: $\alpha=\set{a}$ with $\tau(a)=a$.
In this case $\Pi=\Z/2\Z$ and $
[\Pi,\Pi]=0
$. \par

\subsection{Self-linking of nanowords}\label{sec45}

We introduce here another homotopy  invariant of nanowords, the so-called
self-linking. We begin with the following observation.
Consider the words $ABAB$ and $AABB$. The letters $A,B$ are
obviously linked or interlaced in the first word and unlinked in the second one.
Consider now an arbitrary  nanoword $(\mathcal{A},w)$ over $\alpha$. We say that two
letters $A,B\in\mathcal{A}$ are \emph{$w$--interlaced} if
$$
w=\cdots A\cdots B\cdots A \cdots B\cdots \qquad \text{or} \qquad
w=\cdots B\cdots A\cdots B \cdots A\cdots.
$$
In the first case   set $n_{w}(A,B)=1$ and in the second case
set $n_{w}(A,B)=-1$. In all other cases  set $n_{w}(A,B)=0$. The
function $n_{w}$ is skew-symmetric in the sense that for all $A,B\in\mathcal{A}$,
$$
n_{w}(A,B) =-n_{w}(B,A) \,\,\,\, {\text {and}}\,\,\,\,
n_{w}(A,A) =n_w(B,B)=0.
$$

Now consider the abelian group $
\pi=\Pi\big/[\Pi,\Pi]$. The group operation in $\pi$ will be written
multiplicatively. Each generator $z_a\in \Pi$ with $a\in \alpha$
  projects to an element of $\pi$ denoted
  $a$. Thus,
$$
\pi =\langle\set{a}_{a\in\alpha}\mid ab=ba\text{
and } a\tau(a)=1\text{ for all $a,b\in\alpha$}\rangle.
$$

For a nanoword $(\mathcal{A},w)$ over $\alpha$ and every
$A\in\mathcal{A}$,  set
$$
[A]_{w}=\prod_{B\in\mathcal{A}}\abs{B}^{n_{w}(A,B)}\in\pi.
$$
For   $a\in\alpha$, set
$$
[a]_{w}=\sum_{\ \substack{ A\in\mathcal{A}\text{, }\abs{A}=a\\
[A]_{w}\neq 1}}[A]_{w}\in\Z\pi.
$$
The function $\alpha\mapsto \Z \pi, a\mapsto [a]_w$ is called the \emph{self-linking} of $w$.
The following theorem derives from this function a homotopy
invariant of $w$.

\begin{Thm}\hfill
\begin{enumerate}
\item For any $a\in\alpha$, the difference
$[a]_{w}-[\tau(a)]_{w}\in\Z\pi$ is a homotopy invariant of $w$.
\item If   $\tau(a)=a$, then
$[a]_{w}\bmod{2}\in(\Z\big/2\Z) \pi$ is a homotopy
invariant of~$w$.
\end{enumerate}
\end{Thm}

For a proof, we refer   to \cite{Tu1}.  By this theorem,
the letters of   $\alpha$ give rise to   homotopy invariants of nanowords
over $\alpha$. These
invariants reflect the linking of letters    in   nanowords.
Here is a simple application of this invariant.\par

Pick $a,b\in\alpha$  and consider the nanoword $w=w_{a,b}= ABAB $ with
$\abs{A}=a$ and $\abs{B}=b$. By Lemma \ref{L2}, if   $a=\tau(b)$, then $w  $
is contractible. We can use the invariant $\gamma$ and the self-linking invariant
to show the converse: if    $w $ is   contractible, then $a=\tau(b)$.
Indeed, suppose that $a\neq \tau(b)$.
We have $\gamma(w )=z_{a}z_{b}z_{a}^{-1}z_{b}^{-1}$. If $a\neq
b$, then $a,b$ lie in different orbits of $\tau$ and therefore $z_a$ does not
commute with $z_b$ in $\Pi$. Then  $\gamma(w )\neq 1$
and   $w $ is  non-contractible.
If $a=b\neq \tau(a)$, then $\gamma(w )= 1$. However, in this case
  $\abs{A}=\abs{B}=a=b$ and
\begin{align*}
[A]_{w}&=\abs{A}^{n_{w}(A,A)}\abs{B}^{n_{w}(A,B)}=a, \\
[B]_{w}&=\abs{A}^{n_{w}(B,A)}\abs{B}^{n_{w}(B,B)}=a^{-1}.
\end{align*}
Then  $[a]_{w}=a+a^{-1}\neq
0\in\Z\pi$. Since
there are no letters in $\{A,B\}$ projecting to $\tau_{a}$, we have $[\tau(a)]_{w}=0$.
These computations and the previous theorem imply that $w$ is non-contractible.

\subsection{Applications of the self-linking}

Recall the norm on the nanowords
$$
\norm{w}=\frac{1}{2}(\text{minimal length of a nanoword homotopic to
$w$}).
$$
We can use the self-linking to estimate this norm from below.
The idea  is that elements of a group ring may be treated
like polynomials and for a
polynomial we can consider its degree. When $w$ is not too
long,  there are not so many factors in the self-linking invariant,
and the degree can not be
too big.
Instead of stating here general theorems, we give  an example of the resulting
estimate for a specific nanoword.
Consider the monoliteral  word
$$a^{m}=a  a\cdots\cdot a$$
formed by $m$ copies of a letter  $a\in\alpha$. If $\tau(a)=a$ or $m=1,2$, then this word is
contractible.  If $a\neq\tau(a)$ and $m\geq 3$, then
the
self-linking invariant gives us the following estimate:
$$\norm{a^{m}}=\norm{(a^{m})^{d}}\geq \left[\frac{m}{2}\right]\times
\left[\frac{m-1}{2}\right]+1$$
  where $[x]$ denotes the greatest integer which is smaller than or equal to
$x$. In particular, $a^m$  has a
positive norm and   is non-contractible. The
estimate of $\norm{a^{m}}$ given above is very rough. I
suppose that it gives approximately twice the actual value
of the norm. My conjecture is that
$$
\norm{a^{m}}=\frac{m(m-1)}{2}.
$$

\begin{Thm}
Let $a,b\in\alpha$ such that $a\neq\tau(a)$ and $b\neq\tau(b)$.
  The words $a^{m}$ and $b^{n}$ with $m,n\geq 3$ are homotopic if and
only if $a=b$ and $m=n$.
\end{Thm}

So such   monoliteral  words   are not
homotopic unless they   coincide. The proof goes by
comparing the self-linking invariants.\par

\subsection{Geometric interpretation  of the self-linking}

We give   a geometric interpretation of the self-linking
in the case of curves, that is in the case where   $\alpha$ consists of
two letters $a, b$ permuted by $\tau$. The group $\pi$
is then the infinite cyclic group with generators
$a$, $b$ satisfying $ab=1$.

\begin{figure}[htbp]
\centering \psfrag{A}[cc][cc]{$A_{i}$}\psfrag{i}{}
\includegraphics[scale=1.1]{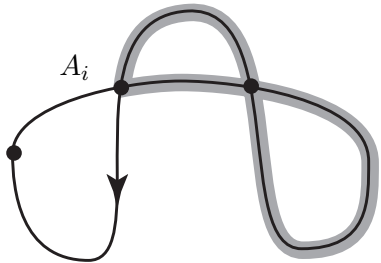}
\caption{Thin and thick curves associated with a crossing $A_i$}
\label{2-28}
\end{figure}

Consider a curve $c$ on an oriented surface with crossings $A_1, A_2, ..., A_n$.   Each crossing $A_i$  gives rise to two sub-curves of $c$  as
follows. Start from   $A_{i}$  and go along $c$ in the positive direction until
coming back to $A_i$ for the first time. The resulting closed curve
is a sub-curve of $c$. The
two branches of $c$ passing through $A_i$ give rise in this way
to two sub-curves
of $c$. One of them passes through the origin of $c$, we
  call this sub-curve the \emph{thin curve}. The other, complementary sub-curve
  is called the \emph{thick curve}, cf.\ Figure \ref{2-28}.

Consider   the  homological intersection number:
$$
k_{i}=\text{thin curve}\cdot\text{thick  curve}\in\Z.
$$
Recall that the homological intersection number of two curves on an oriented surface
is
obtained by     deforming the curves into a transversal
position and then counting their intersections with appropriate signs.

Consider   the nanoword $(\{A_1,..., A_n\}, w)$  corresponding to the curve $c$.
Then for all $i=1,..., n$,
$$
[A_{i}]_{w}=\abs{A_i}^{k_{i}}\in\pi .
$$
This formula gives a geometric interpretation of
   the symbol $[A_i]_w$.  The
self-linking is obtained   by taking the sum   of these symbols over
the crossings
$A_i$ with fixed projection to the alphabet $\{a,b\}$.
To ensure the invariance under
the first Reidemeister move, we   restrict the summation to    $A_i$
such that $k_i\neq 0$ or, equivalently, $[A_i]_w \neq 1$.



\section{Linking   pairings of nanowords}

The invariants   defined so far, $\gamma$ and the self-linking, are   insufficient
  to classify even short words.
We need more invariants.
One idea is to consider again the geometric situation of nanowords associated with
curves. With each crossing $A_i$ we associated a
  `thick curve' on the ambient surface. We can  consider
the intersection numbers  of  these curves with each other. This gives
an $n\times n$ integral matrix
where $n$ is the number of crossings  and the $(i,j)$ entry is
  the intersection number of the thick curves determined  by
$A_{i}$ and  $A_{j}$. This matrix can be computed directly from the nanoword. This
leads us to
  so-called linking pairings of nanowords.  \par

  In this section, the symbol $\pi$ denotes the multiplicative abelian group
  associated with $\alpha, \tau$ in Section \ref{sec45}.

\subsection{$\alpha$--pairings}

We begin with   purely algebraic definitions, whose connection to
nanowords  will be explained later.

\begin{defn}
An \emph{$\alpha$}--pairing is a tuple consisting of  a set $S$, a
distinguished element $s\in S$, a mapping $S-\set{s}\to\alpha$, and
a skew-symmetric pairing $b\co S\times S\to\pi$.
\end{defn}

By skew-symmetric, we mean that   $b(A,B)=b(B,A)^{-1}$ for all $A,B\in S$ and
$b(A,A)=1$ for all $A\in S$.

An $\alpha$--pairing can be shortly written as
$$
(S,s,b\co S\times S\to\pi).
$$
The mapping $S-\set{s}\to\alpha$ will be encoded by saying that the
set $S-\{s\}$ is an $\alpha$--alphabet. The image of any $A\in
S-\set{s}$ under this mapping will be  denoted   $\abs{ A}$.

The notion of isomorphism for  $\alpha$--pairings is defined in the obvious way.

Given an $\alpha$--pairing  $ (S,s,b\co S\times S\to\pi) $, we
define \emph{annihilating elements} of $S$ and \emph{twins} as
follows.

\begin{defn} An element
$A\in S-\set{s}$ is   annihilating    if $b(A,C)=1$ for
all $C\in S$.
\end{defn}

\begin{defn} Elements
$A,B\in S-\set{s}$ are   twins  if $b(A,C)=b(B,C)$ for
all $C\in S$ and $\abs{A}=\tau(\abs{B})$.
\end{defn}

An $\alpha$--pairing is   \emph{primitive} if it has no
annihilating elements and no twins. For example, the \emph{trivial
$\alpha$--pairing} $(S=\{s\},\,
b(s,s)=1)$ is primitive.

We   introduce  two moves  $M_1, M_2$ on $\alpha$--pairings.

\begin{description}
  \item[$M_{1}$] Delete an annihilating element.
  \item[$M_{2}$] Delete a pair of twins.
\end{description}

\noindent The moves $M_{1}$, $M_{2}$,  the  inverse moves $M_{1}^{-1}$, $M_{2}^{-1}$,
  and isomorphisms of
$\alpha$--pairings
generate an equivalence relation    on the class of
$\alpha$--pairings, called \emph{homology}.
The following theorem classifies    $\alpha$--pairings up to homology.

\begin{Thm}\label{pri}
Every $\alpha$--pairing is homologous to a primitive
$\alpha$--pairing. Two homologous primitive $\alpha$--pairings are
isomorphic.
\end{Thm}

Thus,  in each homology class of
$\alpha$--pairings  there is a    primitive one unique up to isomorphism.
Starting with an arbitrary $\alpha$--pairing, we  can
delete annihilating elements and twins and get a primitive
$\alpha$--pairing. The latter  is uniquely determined by the homology class of the
original  $\alpha$--pairing at least up to isomorphism.\par

\subsection{From nanowords to $\alpha$--pairings}\label{fromn}

The connection between $\alpha$--pairings and nanowords  is this:  to
each nanoword $w$ over $\alpha$ we shall
assign an $\alpha$--pairing $b_w$. Its homology class
will be a homotopy
invariant of $w$.

Let
$(\mathcal{A},w\co\widehat {n}\to\mathcal{A})$ be a nanoword over $\alpha$. Set
$S=\set{s}\cup\mathcal{A}$. We have a projection
$S-\set{s}=\mathcal{A}\to\alpha$. The
  skew-symmetric pairing $b_w\co S\times S\to \pi$ is defined in four steps.

\begin{enumerate}[Step 1.]
  \item For every   $A\in\mathcal{A}$, we can write uniquely
  $w^{-1}(A)=\set{i_{A},j_{A}}\subset \widehat  n$ where $i_{A}< j_{A}$. Thus $i_{A}$
   is the position  in which the letter
  $A$   appears in $w$ for the first time, and
  $j_{A}$ is the position  in which the letter $A$  appears in $w$ for the second
   time. Thus $w$ is of the form:
  $$
  w= \cdots \underset{i_{A}}{A}\cdots \underset{j_{A}}{A}\cdots.
  $$
  \item Given two letters $D,E\in\mathcal{A}$, set
  $$D\circ E=
  \prod_{\substack{F\in\mathcal{A}\\ i_{D}<i_{F}<j_{D}\text{ and
  }  i_{E}<j_{F}<j_{E}}}\abs{F}\in\pi.
  $$
  \item The \emph{$w$--linking} of   $D,E\in
  \mathcal{A}$ is defined by
  $$
\mathrm{lk}_{w}(D,E)=(D\circ E)(E\circ D)^{-1}\in\pi.
  $$
  \item Finally,   the form $b_w\co S\times S\to\pi$ is defined  as follows:
  $b_w(s,s)=1$,
$$
b_w(A,s) =[A]_{w}=\prod_{B\in\mathcal{A}}\abs{B}^{n_{w}(A,B)}\in\pi,$$
$$ b_w(s,A) =([A]_{w})^{-1}\in\pi, $$
$$
b_w(A,B) =(\mathrm{lk}_{w}(A,B))^{2}\abs{A}^{n_{w}(A,B)}\abs{B}^{n_{w}(A,B)}\in\pi,$$
for any $A,B\in\mathcal{A}=S-\set{s}$.
\end{enumerate}

The following theorem justifies this definition and
relates  the homotopy of nanowords to the homology of $\alpha$--pairings.

\begin{Thm}
Homotopic nanowords have homologous $\alpha$--pairings.
\end{Thm}

  This theorem together with Theorem \ref{pri},   give  an
efficient method to distinguish
nanowords. Given a nanoword $w$, we first compute  the associated
  $\alpha$--pairing $b_{w}$ and
   then  apply the moves $M_{1}$ and $M_{2}$ to get a
primitive $\alpha$--pairing. The isomorphism class of the latter is a homotopy
invariant of~$w$.

\subsection{Applications}\label{fromn1}
One application of the $\alpha$--pairings is the following estimate of the length norm
of   nanowords: if $(S_+, s, b_+:S_+\times S_+\to \pi)$ is a primitive
$\alpha$--pairing
homologous to the   $\alpha$--pairing $(S, s, b_w)$ of a nanoword $w$,
then $$\Vert w \Vert \geq   {\text {card}} (S_+) -1.$$
Indeed, if $w$ is homotopic to a nanoword $w'$ of length $2m$, then the $\alpha$--pairing
$(S',s, b')$ of $w'$ is homologous to $(S, s, b_w )$ and therefore
reduces by the moves $M_1, M_2$ to
the  same  primitive $\alpha$--pairing $(S_+, s, b_+)$. Hence $$m={\text {card}} (S')
-1\geq {\text {card}} (S_+) -1.$$ In particular, if the $\alpha$--pairing $(S, s, b_w )$
is primitive, then $w$ has minimal length in its homotopy class.

The $\alpha$--pairing $(S,s,b_w)$ can be   used to estimate   the
geometric genera of surfaces carrying the underlying curves of $w$.
Pick an equivariant map $f\co \alpha\to \alpha_0$ where
$\alpha_0=\{a,b\}$ is the 2--letter alphabet with involution
permuting $a,b$. The nanoword $f_{\sharp}(w)$  over $\alpha_0$
(defined in Section \ref{sharrr}) corresponds to a pointed curve on
a compact surface. We can estimate the genus $g$ of this surface
  by
$g \geq (1/2)\, {\text {rank}} (M)$ where $M$ is the  skew-symmetric integral
matrix obtained from the
matrix $\{b_w(s_1, s_2)\}_{s_1, s_2 \in S}$ by the group homomorphism $\pi\to \Z$
sending   the generators of $\pi$ belonging to $f^{-1}(a)$ to $1\in \Z$ and the
generators of $\pi$ belonging to $f^{-1}(b)$ to $-1\in \Z$. This estimate follows from
the geometric interpretation of the $\alpha_0$-pairing  of $f_{\sharp}(w)$ in terms of  the
intersection numbers of curves.

Another area of applications of  $\alpha$--pairings is the homotopy
classification of nanowords.
With the help of $\alpha$--pairings we can establish the following theorem.
Recall the nanoword $w_{a,b}= ABAB,
\abs{A}=a, \abs{B}=b $ defined for any $a,b\in \alpha$. As we know, $w_{a,b}$ is
non-contractible if and only if $a\neq\tau(b)$.

\begin{Thm}
Two
non-contractible  nanowords $w_{a,b}$ and $w_{a',b'}$  with $a,b,a'$, $b'\in \alpha$
are
homotopic if and only if $a=a'$ and $b=b'$.
\end{Thm}

Using the $\alpha$--pairings and the invariants introduced in further sections, we
establish the following theorem. It gives  a complete homotopy
classification of      words of length $5$ in which one letter, $a$, occurs 3 times,
and another letter, $b$, occurs 2 times.

\begin{Thm}\label{T:len5word}
Let $a,b$ be two distinct letters of the alphabet $\alpha$. Then:
\begin{enumerate}
\item The   words $aaabb$, $aabba$, $abbaa$, $bbaaa$ are homotopic to each other; they
are contractible if and only if $\tau (a)=a$.
\item The word $baaab$ is contractible if and only if $\tau(a)=a$.
\item The word $ababa$ is contractible if and only if $\tau(a)=b$.
\item The   words $abaab$, $baaba$, $aabab$,
$babaa$  are never contractible.
\item A non-contractible word from   (2) -- (4) is never homotopic to a word from (1).
\item Two non-contractible words from  (2) -- (4)   are homotopic to each other  if
and only if they coincide letterwise (i.e.,  if and only if
they are the same word written twice) with the following exceptions:
$$
aabab\simeq babaa\simeq baaab \qquad\qquad\text{for $ \tau(a)=b$.}
$$
\end{enumerate}
\end{Thm}

  A more general homotopy  classification of
all words of length $\leq 5$ is given in \cite{Tu1}.
We can think of such classification theorems as analogues of knot tables.
First we draw all  possible knot diagrams  and then decide
which diagrams   represent isotopic knots. The
same kind of problem arises for the homotopy of words.\par

\subsection{Examples}

1. We show how to compute the $\alpha$--pairing associated with the
word
  $w=abaab$ where $a\neq b$. We have
$$
w^{d}= A_{3}A_{2}BA_{3}A_{1}A_{2}A_{1}B,\,\,\,
\abs{A_{1}}=\abs{A_{2}}=\abs{A_{3}}=a,\ \abs{B}=b
$$
  where to simplify notation we write $A_{1},A_{2},
  A_{3}$ for  $A_{2,3}, A_{1,3}, A_{1,2}$
respectively. The matrix for $n_{w}$ is computed by
\begin{equation}
\left[
\begin{array}{cccc}
\phantom{-}0          & -1            & \phantom{-}0            & \phantom{-}0 \\
\phantom{-}1          & \phantom{-}0            & -1           & \phantom{-}1\\
\phantom{-}0          & \phantom{-}1            & \phantom{-}0            & \phantom{-}1\\
\phantom{-} 0         & -1            & -1             &
\phantom{-}0
\end{array}
\right]
\notag
\end{equation}
where the rows and columns correspond to $A_1, A_2, A_3,
B$ respectively. The matrix for $\mathrm{lk}_{w}$ is computed by
\begin{equation}
\left[
\begin{array}{llll}
1          & 1            & a^{-1}            & 1 \\
1          & 1            & 1           & a\\
a          & 1            & 1            & a\\
1        & a^{-1}           & a^{-1}             & 1
\end{array}
\right]
\notag
\end{equation}
  where the rows and columns correspond to $A_1, A_2, A_3,
B$, respectively. Finally, we compute the $\alpha$--pairing $b_w$:
\begin{equation}
b_w=\left[
\begin{array}{lllll}
1  \phantom{a^{-1}}           & a \phantom{a^{-1}}         & b^{-1}  \phantom{a^{-1}}           & a^{-1}b^{-1} & a^{2} \phantom{a^{-1}}\\
a^{-1} \phantom{a^{-1}} & 1 \phantom{a^{-1}}          & a^{-2} & a^{-2} & 1 \\
b \phantom{a^{-1}}            & a^{2} \phantom{a^{-1}} & 1 \phantom{a^{-1}}            & a^{-2} \phantom{a^{-1}} & a^{3}b\\
ab           & a^2  \phantom{a^{-1}}    & a^2  \phantom{a^{-1}}      & 1 \phantom{a^{-1}}            & a^{3}b\\
a^{-2} \phantom{a^{-1}}            & 1 \phantom{a^{-1}}         & a^{-3}b^{-1}            & a^{-3}b^{-1} & 1
\end{array}
\right]
\notag
\end{equation}
where the rows and columns correspond to $s, A_1, A_2,
A_3, B$, respectively. Recall that $\pi=\langle \{c\}_{c \in
\alpha}\mid c\tau(c)=1 \rangle$.
In particular $a$ and $b$ are non-trivial elements of
$\pi$. This implies that the elements $A_{1}$ and  $A_{2}$  of
$S=\{s, A_1, A_2, A_3,
B\}$
are   non-annihilating. If $A_{3}$ is annihilating, then
$a^{2}=ab=1$. Then $ a=b$, which contradicts the assumption $a\neq b$.
A similar argument shows that $B$ is non-annihilating. It is
also easy to check that $b_w$ does not have  twins. Thus,
the $\alpha$--pairing $b_w$ is primitive. Therefore $w$ is
non-contractible, and $\norm{w}=4$.\par

2. The $\alpha$--pairings are   strong enough
to distinguish short  words and nanowords in many cases. The following example shows
however that in some cases the $\alpha$--pairings are powerless.

Consider  the  word $w=ababa$ where
$ \tau(a)=a \neq b=  \tau(b) $.
A direct computation shows that the $\alpha$--pairing of $w^d$ is
  given by the following matrix over   $\pi$:
\begin{equation}
\left[
\begin{array}{ccccc}
1            & ab         & 1            & a^{-1}b^{-1} & 1\\
a^{-1}b^{-1} & 1          & a^{-2}b^{-2} & a^{-2}b^{-2} & a^{-1}b^{-1} \\
1            & a^{2}b^{2} & 1            & a^{-2}b^{-2} & 1\\
ab           & a^2b^2     & a^2b^2       & 1            & ab\\
1            & ab         & 1            & a^{-1}b^{-1} & 1
\end{array}
\right].
\notag
\end{equation}
\noindent As above, the rows and columns correspond to $s, A_1, A_2, A_3, B$
respectively.    The equality $a=\tau(a)$ implies that $a^2=1$  and
similarly  $b^2=1$. Therefore the matrix above simplifies to the following matrix:
\begin{equation}
\left[
\begin{array}{ccccc}
1            & ab         & 1            & a b  & 1\\
ab & 1          & 1            & 1            & ab \\
1            & 1          & 1            & 1            & 1\\
ab           & 1          & 1            & 1            & ab\\
1            & ab         & 1            & a b  & 1
\end{array}
\right].
\notag
\end{equation}
Since the third row  and the third column  consist  only of   $1's$, the element
$A_2$ is   annihilating.
Eliminating it, we observe next that
$A_1$ and $A_3$ are twins. What
remains after their elimination are two elements $s$ and $B$. The matrix becomes as
follows:
\begin{equation} \left[
\begin{array}{cc}
1 & 1\\
1 & 1
\end{array}
\right]
\notag
\end{equation}
where the rows and columns correspond to $s$ and $B$.
Now $B$ is an annihilating element. Its elimination gives the trivial
$\alpha$--pairing.   Therefore the
$\alpha$--pairing associated with $w^d$
gives no information at all and does not allow to decide whether
$w=ababa$ is contractible or not. In fact all invariants of nanowords
considered so far
are trivial for this word (under the assumptions that
$ \tau(a)=a \neq b=  \tau(b) $).

This example shows that   we need more
invariants  to prove
Theorem \ref{T:len5word}. We shall introduce further invariants in the next sections.

\section{Further invariants of nanowords}\label{fgi134}

We  outline here several ideas  inspired by knot theory and leading to homotopy invariants
of nanowords over $\alpha$.

\subsection{Tricolorings}\label{tri} In knot theory one may treat any knot
diagram as consisting of disjoint  arcs. A \emph{coloring} of the
diagram assigns a residue mod 3
to each arc, such that for every crossing, the sum of the three residues assigned
to the adjacent arcs   is equal to $0$. The  number of such colorings
of a knot diagram is a knot invariant. This definition is due to R.\ Fox.

We can introduce similar definitions  for words.
Fix a set $\beta \subset \alpha$ such that  $\tau(\beta)=\beta$  (the resulting invariants may depend on
$\beta$). Consider a nanoword $w=(\mathcal A,
w:\widehat {n}\to\mathcal A)$ over $\alpha$. For any letter $A\in \mathcal A$, let $i_A< j_A$ be the first and the second
indices   enumerating the positions of
$A$ in   $w$ as in Section \ref{fromn}.
A {\it tricoloring} of $w$ is a function
$f\co\{0,1,2,\ldots,n\} \to \Z/3\Z$ such that for any
$A\in\mathcal A$, if $|A|\in\beta$, then
$$
f(i_A)=f(i_A-1) \quad {\text {and}} \quad
f(j_A-1)+f(j_A)+f(i_A)=0$$
and if $|A|\in\alpha-\beta$, then
$$
f(j_A)=f(j_A-1) \quad {\text {and}} \quad
f(i_A-1)+f(i_A)+f(j_A)=0.$$
The residues $f(0)$ and $f(n)$ are called respectively the \emph{input} and
the \emph{output} of  $f$.

Tricolorings of $w$ may be alternatively described as follows. We first write $w$
with dashes between consecutive letters:
\begin{equation}
- \, w(1)-w(2)-\cdots-w(n)\, -
\notag
\end{equation}
Enumerate the dashes from left to right  by the numbers $0,1,
\ldots, n$. Then the function $f$ as above can be seen as an
assignment of a residue  mod 3 to every dash.
  The
  conditions above mean  that for any
$A\in\mathcal A$, the coloring  has the following form
near the two entries of $A$ in $w$:
if $|A|\in\beta$, then  it looks like
\[ \cdots\dash{c}A\dash{c}\cdots\dash{c'}A\dash{c''}\cdots
\]
with $c+c'+c''=0$ and if  $|A|\in \beta- \alpha$, then it looks like
\[
\cdots\dash{c'}A\dash{c''}\cdots\dash{c}A\dash{c}\cdots
\]
with $c+c'+c''=0$. The input of the coloring is the residue
  assigned to the leftmost dash and the output is the residue   assigned to the
   rightmost dash.

For example, the function assigning one and the same residue to all dashes is a coloring.
It is called the \emph{trivial coloring}.

\begin{Thm}\label{Thm:tricoloring}
For any $k,l\in \Z/3\Z$, the number of tricolorings  of a nanoword $w$ with   input  $k$
and   output $l$  is a homotopy invariant of $w$.
\end{Thm}

Note that this number may depend  on $k, l$, and on the choice of $\beta$.
We can easily compute this number
for the empty nanoword (it has only   one dash).
If $k=l$, then this nanoword admits one coloring with input $k$ and output $l$. If
$k\neq l$, then there are no such colorings.
By Theorem \ref{Thm:tricoloring}, the same is true for any contractible nanoword.

Now we give an example of a nanoword which admits a coloring with distinct input and
output. Consider the nanoword
\begin{equation}
w=A_1A_2BA_3A_1BA_2A_3,
\notag
\end{equation}
where $|A_1|=|A_2|=|A_3|=a\in \alpha$,  $|B|=b\in \alpha$
such that $a$ and $b$ lie in different orbits of $\tau$. Set $\beta=\{a,\tau(a)\}\subset
\alpha$. The nanoword $w$ has the following non-trivial coloring with input 0 and output 1:
\[
\dash{0}A_1\dash{0}A_2\dash{0}B\dash{1}A_3\dash{1}A_1\dash{2}B\dash{2}A_2
\dash{1}A_3\dash{1}.
\]
By the remarks above, this
nanoword is non-contractible.

\subsection{Module of a nanoword}\label{psiintro}
A related but stronger invariant of knots is the Alexander module.
It can be computed from a  knot diagram
via an explicit presentation by generators and relations.
We   apply   a similar idea to  nanowords.
First, we introduce  the group
\[ \Psi=\langle \{a,a_\bullet\}_{a\in\alpha}\,|\,aa_\bullet=a_\bullet a, \,
a\tau(a)=1, \,a_\bullet\tau(a)_\bullet=1\rangle. \]
We   already considered groups $\Pi$ and $\pi$ given by similar presentations. In $\Psi$,
each
letter $a\in \alpha$ gives rise to two commuting generators $a$ and $a_\bullet$.
  In the case of knot diagrams on surfaces, this phenomenon of doubling of the number of
  generators was
already observed by Sawollek \cite{SA} who studied generalizations of the Alexander
polynomial, see also   \cite{SW1}.

Let $\Lambda=\Z\Psi$ be the integral group ring of
$\Psi$. This ring will play  the role of the ground ring for our modules.

Fix a set $\beta\subset\alpha$ such that $\tau(\beta)=\beta$.
Consider a nanoword $w=(\mathcal A, w\co\widehat {n}\to\mathcal A)$
over $\alpha$.
  We derive from $w$  a $(n+1)\times n$ matrix over $\Lambda$ whose rows are  numerated
by the dashes of $w$.  Each
letter $A\in\mathcal A$  gives rise to two rows. To write them down,
set $a=|A|\in\alpha$ and assume that $A$ appears in $w$ for the first time
at the $i$-th position and   for the second time at the $j$-th
position.
If $a\in\beta$, then the two rows determined by $A$ are
\begin{equation}
\begin{array}{ccccc}
&\scriptstyle i&&\scriptstyle j&\\
\cdots & -\ \ A\ \ - & \cdots & -\ \ A\ \ - & \cdots\\
    \begin{array}{c}  \cdots \,  0\\  \cdots \, 0\end{array} &
    {\begin{array}{cc}
        a & -1\\
        1-aa_\bullet & 0
        \end{array}} &
    {\begin{array}{c}0\, \cdots \, 0\\0\, \cdots \, 0\end{array}} &
    {\begin{array}{cc}
        0 & 0\\
        a_\bullet & -1
        \end{array}} &
     {\begin{array}{c}0\, \cdots  \\0\, \cdots   \end{array}}
\end{array}
\notag
\end{equation}

If $a\in\alpha-\beta$, then the two rows determined by $A$ are
\begin{equation}
\begin{array}{ccccc}
&\scriptstyle i&&\scriptstyle j&\\
\cdots & -\ \ A\ \ - & \cdots & -\ \ A\ \ - & \cdots\\
    \begin{array}{c}\cdots \,  0\\  \cdots \, 0\end{array} &
    {\begin{array}{cc}
        0 & 0\\
        a_\bullet & -1
        \end{array}} &
    {\begin{array}{c}0\, \cdots \, 0\\0\, \cdots \, 0\end{array}} &
    {\begin{array}{cc}
        a & -1\\
        1-aa_\bullet & 0
        \end{array}} &
     {\begin{array}{c}0\, \cdots  \\0\, \cdots  \end{array}}
\end{array}
\notag
\end{equation}
All unspecified entries of the rows are 0.
The resulting
$n\times (n+1)$ matrix over $\Lambda$
determines a $\Lambda$-homomorphism $\psi\co\Lambda^n \to \Lambda^{n+1}$ whose cokernel
\[ K_{\beta}(w)=\Lambda^{n+1}/\psi(\Lambda^n). \]
\noindent is a $\Lambda$--module. This module has
distinguished elements: the ``input'' $ v_- $  and the ``output''
$v_+ $  represented by the leftmost dash and the rightmost dash, respectively.

\begin{Thm}\label{Thm:module_invariant}
The triple $(K_{\beta}(w),v_-,v_+)$ considered up to isomorphism is
a homotopy invariant of $w$.
\end{Thm}

One can  derive further homotopy invariants from  the triple  $(K_{\beta}(w),v_-,v_+)$
or directly from the presentation matrix of $K_{\beta}(w)$ introduced above.
For example, one may remove the first (or the last) column and consider the
resulting $n\times n$
matrix over $\Lambda$. Then we can take its determinant over the ring
$\Lambda^{ab}$ obtained by abelianization of
$\Lambda$. Following this line of thought  and with a little more work, one obtains two
\lq\lq polynomial" homotopy invariants of $w$ belonging to $\Lambda^{ab}$ (see
\cite{Tu1} for
details). They are denoted
$\nabla^-_{\beta}(w)$  and $\nabla^+_{\beta}(w)$ and satisfy the following   duality:
\begin{equation}
\nabla^+_{\beta}(w)=\overline{\nabla^-_{\alpha-\beta}(w^-)}. \notag
\end{equation}
The bar on the right-hand side is the ring involution on  $\Lambda^{ab}$ given by
$
\overline{a}=\tau(a)$ and $\overline{a_\bullet}=\tau(a)_\bullet
$ for all $a\in \alpha$.

\subsection{The invariant $\lambda$}
We now focus on the case where $\beta=\alpha$. This will  lead us to a homotopy
  invariant
$\lambda$ of nanowords taking values in the ring $\Lambda$. This invariant is a
generalization of an invariant
introduced by Silver and Williams \cite{SW2} for curves.

Consider the $\Lambda$--module $K_{\beta}(w)=K_{\alpha}(w)$
associated   above with a nanoword $w\co\widehat {n}\to\mathcal A$.
Let $x_0,x_1,\ldots,x_n$ be the generators of
  $K_{\alpha}(w)$
given by the dashes. Each letter  $A\in\mathcal A$
gives rise to two relations
$$
x_{i+1} =ax_{i},\quad
x_{j} =a_\bullet x_{j-1}+(1-aa_\bullet)x_{i},
$$
where $a= |A| \in\alpha$ and $i=i_A < j=j_A$  are the  elements of $w^{-1}(A)$.
Each of these relations expresses a generator via the previous
generators. Therefore   $K_{\alpha}(w)$ is   a
rank one free $\Lambda$-module generated by the input $v_-=x_0$.
The  output $v_+=x_n\in K_{\alpha}(w)$ has
the form $v_+=\lambda' v_- $ for a unique $\lambda'\in \Lambda$.
Theorem \ref{Thm:module_invariant}
implies    that $\lambda'=\lambda'(w)$ is a
homotopy invariant of $w$. This invariant is a non-commutative
polynomial. It admits an equivalent but more convenient version $\lambda(w)$
defined as follows. Consider the involutive anti-automorphism $\iota$ of $\Lambda$
keeping fixed all the   generators $\{a, a_\bullet\}$ of $\Lambda$.
Thus, $\iota$ acts on
monomials by reading them from right to left. For instance $\iota(aab_\bullet)=
b_\bullet aa$. Set
$$\lambda(w)=\iota (\lambda'(w))\in \Lambda.$$

We describe a method allowing to compute $\lambda(w)$ and generalizing a method due
to Silver and Williams \cite{SW2} in the context of  curves.
We   do it here  for a few
examples, the general method   \cite{Tu1} should be clear.

\begin{Exa} Consider the nanoword $w=ABAB,\ |A|=a\in\alpha,\ |B|=b\in\alpha$.
First draw the following graph:

\begin{center}
\unitlength 0.1in
\begin{picture}( 24.8000,  5.1000)(  5.3000,-12.8000)
\put(14.3000,-9.4000){\makebox(0,0)[lb]{$b$}}%
%
\special{pn 8}%
\special{pa 610 1000}%
\special{pa 3010 1000}%
\special{fp}%
\put(5.7000,-12.0000){\makebox(0,0)[lb]{0}}%
\put(11.6000,-12.0000){\makebox(0,0)[lb]{1}}%
\put(17.6000,-12.0000){\makebox(0,0)[lb]{2}}%
\put(23.5000,-12.0000){\makebox(0,0)[lb]{3}}%
\put(29.5000,-12.0000){\makebox(0,0)[lb]{4}}%
\put(8.3000,-9.4000){\makebox(0,0)[lb]{$a$}}%
\put(20.0000,-9.4000){\makebox(0,0)[lb]{$a_\bullet$}}%
\put(26.0000,-9.4000){\makebox(0,0)[lb]{$b_\bullet$}}%
%
\special{pn 8}%
\special{sh 1}%
\special{ar 600 1000 10 10 0  6.28318530717959E+0000}%
\special{sh 1}%
\special{ar 600 1000 10 10 0  6.28318530717959E+0000}%
%
\special{pn 8}%
\special{sh 1}%
\special{ar 1200 1000 10 10 0  6.28318530717959E+0000}%
\special{sh 1}%
\special{ar 1200 1000 10 10 0  6.28318530717959E+0000}%
%
\special{pn 8}%
\special{sh 1}%
\special{ar 1800 1000 10 10 0  6.28318530717959E+0000}%
\special{sh 1}%
\special{ar 1800 1000 10 10 0  6.28318530717959E+0000}%
%
\special{pn 8}%
\special{sh 1}%
\special{ar 2400 1000 10 10 0  6.28318530717959E+0000}%
\special{sh 1}%
\special{ar 2400 1000 10 10 0  6.28318530717959E+0000}%
%
\special{pn 8}%
\special{sh 1}%
\special{ar 3000 1000 10 10 0  6.28318530717959E+0000}%
\special{sh 1}%
\special{ar 3000 1000 10 10 0  6.28318530717959E+0000}%
\put(5.3000,-14.5000){\makebox(0,0)[lb]{---}}%
\put(10.9000,-14.5000){\makebox(0,0)[lb]{---}}%
\put(16.9000,-14.5000){\makebox(0,0)[lb]{---}}%
\put(23.0000,-14.4000){\makebox(0,0)[lb]{---}}%
\put(29.2000,-14.5000){\makebox(0,0)[lb]{---}}%
\put(8.4000,-14.3000){\makebox(0,0)[lb]{$A$}}%
\put(14.5000,-14.3000){\makebox(0,0)[lb]{$B$}}%
\put(20.2000,-14.3000){\makebox(0,0)[lb]{$A$}}%
\put(26.5000,-14.2000){\makebox(0,0)[lb]{$B$}}%
\end{picture}%
\vspace{7mm}
\end{center}
Each vertex of this graph corresponds to a dash in $w$ and each edge corresponds to a
letter in $w$. Recall that every letter appears twice. The edge corresponding to the first
(leftmost)
appearance of $A$ is labeled with $a$; the edge corresponding to the second
(rightmost)
appearance of $A$ is labeled with $a_\bullet$. Connect the
left vertex of the first edge with the right vertex of the second edge by an arc in the
upper half-plane and label this arc with $1-aa_\bullet\in \Lambda$.
Do the same for the letter $B$ replacing everywhere $a=\abs{A}$ by $b=\abs{B}$.  The
resulting picture is drawn on the next figure.

\begin{center}
\unitlength 0.1in
\begin{picture}( 28.8000,  8.3000)( 14.0000,-24.3000)
\put(26.3000,-23.4000){\makebox(0,0)[lb]{$b$}}%
%
\special{pn 8}%
\special{pa 1810 2400}%
\special{pa 4210 2400}%
\special{fp}%
\put(17.7000,-26.0000){\makebox(0,0)[lb]{0}}%
\put(23.6000,-26.0000){\makebox(0,0)[lb]{1}}%
\put(29.6000,-26.0000){\makebox(0,0)[lb]{2}}%
\put(35.5000,-26.0000){\makebox(0,0)[lb]{3}}%
\put(41.5000,-26.0000){\makebox(0,0)[lb]{4}}%
\put(20.3000,-23.4000){\makebox(0,0)[lb]{$a$}}%
\put(32.0000,-23.4000){\makebox(0,0)[lb]{$a_\bullet$}}%
\put(38.0000,-23.4000){\makebox(0,0)[lb]{$b_\bullet$}}%
%
\special{pn 8}%
\special{sh 1}%
\special{ar 1800 2400 10 10 0  6.28318530717959E+0000}%
\special{sh 1}%
\special{ar 1800 2400 10 10 0  6.28318530717959E+0000}%
%
\special{pn 8}%
\special{sh 1}%
\special{ar 2400 2400 10 10 0  6.28318530717959E+0000}%
\special{sh 1}%
\special{ar 2400 2400 10 10 0  6.28318530717959E+0000}%
%
\special{pn 8}%
\special{sh 1}%
\special{ar 3000 2400 10 10 0  6.28318530717959E+0000}%
\special{sh 1}%
\special{ar 3000 2400 10 10 0  6.28318530717959E+0000}%
%
\special{pn 8}%
\special{sh 1}%
\special{ar 3600 2400 10 10 0  6.28318530717959E+0000}%
\special{sh 1}%
\special{ar 3600 2400 10 10 0  6.28318530717959E+0000}%
%
\special{pn 8}%
\special{sh 1}%
\special{ar 4200 2400 10 10 0  6.28318530717959E+0000}%
\special{sh 1}%
\special{ar 4200 2400 10 10 0  6.28318530717959E+0000}%
%
\special{pn 8}%
\special{ar 2700 2400 900 800  3.1415927 6.2831853}%
%
\special{pn 8}%
\special{ar 3300 2400 900 800  3.1415927 6.2831853}%
\put(13.8000,-19.2000){\makebox(0,0)[lb]{$1-aa_\bullet$}}%
\put(40.4000,-19.2000){\makebox(0,0)[lb]{$1-bb_\bullet$}}%
\put(14.0000,-24.7000){\makebox(0,0)[lb]{input}}%
\put(42.8000,-24.7000){\makebox(0,0)[lb]{output}}%
\end{picture}%
\vspace{7mm}
\end{center}

Consider all paths starting at the input and going to the
output along the  edges and arcs, always
  from left to right. We record the elements of $\Lambda$
labeling the  arcs and edges on the path  and multiply  them following the
order determined by the path. The polynomial $\lambda(w)$ is obtained as the sum of
the resulting elements of $\Lambda$ over all paths.
In this case there are three such paths:
\begin{enumerate}
\item The path $0-1-2-3-4$  contributes $aba_\bullet b_\bullet$.
\item The path  $0-1\frown 4$ contributes  $a(1-bb_\bullet)$.
\item The path  $0\frown 3-4$  contributes $(1-aa_\bullet)b$.
\end{enumerate}

\noindent Then
\begin{equation}
\lambda(w)=aba_\bullet b_\bullet+a(1-bb_\bullet)+(1-aa_\bullet)b_\bullet.
\notag
\end{equation}

The ring $\Lambda$   has a natural
grading as follows.
Recall that
\begin{equation}
\Lambda=\Z[a,a_\bullet]_{a\in\alpha}/
aa_\bullet=a_\bullet a,\ a\tau(a)=1,\ a_\bullet\tau(a)_\bullet=1 .
\notag
\end{equation}
The  defining relations are homogeneous with respect
to degrees mod 2. Therefore
\begin{equation}
\Lambda=\Lambda_{0,0}\oplus\Lambda_{0,1}\oplus\Lambda_{1,0}\oplus\Lambda_{1,1},
\notag
\end{equation}
where $\Lambda_{i,j}$ is generated by monomials in which generators
without bullets appear $i$ times mod  $2$ and generators with bullets appear
$j$ times mod  $2$. Every   $\lambda\in\Lambda$ expands uniquely as
the   sum
\begin{equation}
\lambda=\lambda_{0,0}+\lambda_{0,1}+\lambda_{1,0}+\lambda_{1,1},
\notag
\end{equation}
where $\lambda_{i,j}\in \Lambda_{i,j}$ for all $i,j$.
For     $w=ABAB$, this expansion of $\lambda(w)$ gives:
\begin{equation}
\begin{split}
\lambda_{0,0}(w)&=aba_\bullet b_\bullet,\\
\lambda_{0,1}(w)&=-abb_\bullet+b_\bullet,\\
\lambda_{1,0}(w)&=a-aa_\bullet b_\bullet,\\
\lambda_{1,1}(w)&=0.
\end{split}
\notag
\end{equation}
These computations allow us to give another  proof of the fact that $w$ is contractible if  and only if
$a=\tau(b)$. Indeed, if  $w$ is contractible,  then $\lambda (w)=1$ and hence
$\lambda_{0,1}(w)=
0$.
This implies  that
$abb_\bullet=b_\bullet$. Hence $ab=1$ and
$a=\tau(b)$.
\end{Exa}

\begin{Exa}
We   apply  $\lambda$ to the word $ababa$ with $a=\tau(a)\neq b=\tau(b)$. As we saw above, the
corresponding $\alpha$--pairing gives no information about the homotopy properties of $w$.
By definition, $\lambda(w)=\lambda(w^d)$.  The desingularization of $w$ is the nanoword
\begin{equation} w^d=A_3A_2BA_3A_1BA_2A_1,\
|A_1|=|A_2|=|A_3|=a,\ |B|=b. \notag
\end{equation}
To
compute   $\lambda(w)$, we   draw the following graph:\\
\begin{center}
\unitlength 0.1in
\begin{picture}( 33.6600, 12.8600)( 14.5000,-23.1000)
\put(25.8700,-22.4300){\makebox(0,0)[lb]{$b$}}%
%
\special{pn 8}%
\special{pa 1680 2270}%
\special{pa 3248 2270}%
\special{fp}%
\put(18.1600,-22.4300){\makebox(0,0)[lb]{$a$}}%
\put(41.6200,-22.4300){\makebox(0,0)[lb]{$a_\bullet$}}%
\put(37.7000,-22.4300){\makebox(0,0)[lb]{$b_\bullet$}}%
%
\special{pn 8}%
\special{sh 1}%
\special{ar 1680 2270 10 10 0  6.28318530717959E+0000}%
\special{sh 1}%
\special{ar 1680 2270 10 10 0  6.28318530717959E+0000}%
%
\special{pn 8}%
\special{sh 1}%
\special{ar 2072 2270 10 10 0  6.28318530717959E+0000}%
\special{sh 1}%
\special{ar 2072 2270 10 10 0  6.28318530717959E+0000}%
%
\special{pn 8}%
\special{sh 1}%
\special{ar 2464 2270 10 10 0  6.28318530717959E+0000}%
\special{sh 1}%
\special{ar 2464 2270 10 10 0  6.28318530717959E+0000}%
%
\special{pn 8}%
\special{sh 1}%
\special{ar 2856 2270 10 10 0  6.28318530717959E+0000}%
\special{sh 1}%
\special{ar 2856 2270 10 10 0  6.28318530717959E+0000}%
%
\special{pn 8}%
\special{sh 1}%
\special{ar 3248 2270 10 10 0  6.28318530717959E+0000}%
\special{sh 1}%
\special{ar 3248 2270 10 10 0  6.28318530717959E+0000}%
\put(13.5000,-16.8200){\makebox(0,0)[lb]{$1-aa_\bullet$}}%
\put(30.1900,-15.4600){\makebox(0,0)[lb]{$1-bb_\bullet$}}%
%
\special{pn 8}%
\special{sh 1}%
\special{ar 3640 2270 10 10 0  6.28318530717959E+0000}%
\special{sh 1}%
\special{ar 3640 2270 10 10 0  6.28318530717959E+0000}%
%
\special{pn 8}%
\special{sh 1}%
\special{ar 4032 2270 10 10 0  6.28318530717959E+0000}%
\special{sh 1}%
\special{ar 4032 2270 10 10 0  6.28318530717959E+0000}%
%
\special{pn 8}%
\special{sh 1}%
\special{ar 4424 2270 10 10 0  6.28318530717959E+0000}%
\special{sh 1}%
\special{ar 4424 2270 10 10 0  6.28318530717959E+0000}%
%
\special{pn 8}%
\special{sh 1}%
\special{ar 4816 2270 10 10 0  6.28318530717959E+0000}%
\special{sh 1}%
\special{ar 4816 2270 10 10 0  6.28318530717959E+0000}%
%
\special{pn 8}%
\special{pa 3248 2270}%
\special{pa 4816 2270}%
\special{fp}%
\put(22.0200,-22.4300){\makebox(0,0)[lb]{$a$}}%
\put(29.8600,-22.4300){\makebox(0,0)[lb]{$a_\bullet$}}%
\put(33.7800,-22.4300){\makebox(0,0)[lb]{$a$}}%
\put(45.5400,-22.4300){\makebox(0,0)[lb]{$a_\bullet$}}%
%
\special{pn 8}%
\special{ar 2466 2270 788 724  3.1415927 6.2831853}%
%
\special{pn 8}%
\special{ar 4032 2270 784 724  3.1415927 6.2831853}%
%
\special{pn 8}%
\special{ar 3248 2270 786 724  3.1415927 6.2831853}%
%
\special{pn 8}%
\special{ar 3248 2270 1174 1068  3.1415927 6.2831853}%
\put(29.9200,-11.9400){\makebox(0,0)[lb]{$1-aa_\bullet$}}%
\put(45.2800,-16.8200){\makebox(0,0)[lb]{$1-aa_\bullet$}}%
\put(25.8100,-24.8000){\makebox(0,0)[lb]{$B$}}%
\put(18.1000,-24.8000){\makebox(0,0)[lb]{$A_3$}}%
\put(41.5600,-24.8000){\makebox(0,0)[lb]{$A_2$}}%
\put(37.6400,-24.8000){\makebox(0,0)[lb]{$B$}}%
\put(21.9600,-24.8000){\makebox(0,0)[lb]{$A_2$}}%
\put(29.8000,-24.8000){\makebox(0,0)[lb]{$A_3$}}%
\put(33.7200,-24.8000){\makebox(0,0)[lb]{$A_1$}}%
\put(45.4800,-24.8000){\makebox(0,0)[lb]{$A_1$}}%
\end{picture}%
\vspace{7mm}
\end{center}
Then
\begin{equation}
\begin{array}{lll}
    \lambda(w)\,=& (1-aa_\bullet)^2 & \frown\frown\\
        & +(1-aa_\bullet)ab_\bullet a_\bullet^2 & \frown----\\
        & +a(1-aa_\bullet)a_\bullet & -\frown-\\
        & +a^2(1-bb_\bullet)a_\bullet^2 & --\frown--\\
        & +a^2ba_\bullet(1-aa_\bullet) & ----\frown\\
        & +a^2baa_\bullet b_\bullet a_\bullet^2. & --------
\end{array}
\notag
\end{equation}
The assumptions $a=\tau(a)$ and $b=\tau(b)$
imply that $a^2= b^2=a^2_\bullet=b^2_\bullet=1$.  After
  simplification, we obtain that
\begin{equation}
\lambda_{0,0}(w)=2-ba-a_\bullet b_\bullet+baa_\bullet b_\bullet.
\notag
\end{equation}
If $\lambda_{0,0}(w)=1$, then one of the two elements $ba$ and
$a_\bullet b_\bullet$ of the group $\Psi$ must  be equal to $1$. This is possible
only if $a=\tau(b)=b$, which contradicts the assumptions. So,  $\lambda_{0,0}(w)
\neq 1$  and
  $w$ is non-contractible.
\end{Exa}

\begin{Exa}
Consider the
nanowords
\begin{equation}
\begin{split}
w_1&=ABACBC,\ |A|=|C|=a,\ |B|=  \tau(a),\\
w_2&=ACAC,\ |A|=|C|=a,
\end{split}
\notag
\end{equation}
where $a\in \alpha$ satisfies $\tau(a)\neq a$. These two nanowords
are not  distinguished by $\lambda$. In fact,   all the techniques described so far fail to distinguish these
nanowords up to homotopy. This   can be done using the methods introduced
  in the next section. \end{Exa}
\section{$\alpha$--keis and     words}
\subsection{$\alpha$--keis}
Keis were introduced in 1942 by a Japanese mathematician, M.\ Takasaki,   see S.\
Kamada \cite{Kam}
for a comprehensive survey of
keis, their generalizations, and connections with knot theory.
A {\it kei} is a set $X$ with   multiplication $*$ which   satisfies a few  axioms, the main axiom being
\begin{equation}
(x*y)*z=(x*z)*(y*z)  \notag
\end{equation}
for all $x,y,z\in X$. One may think of $x*y$ as  of a kind of conjugation of $x$ by $y$.

  To produce homotopy invariants of words, we introduce a notion of an $\alpha$--kei,
   where   $\alpha$ is a set with involution $\tau$.
An {\it $\alpha$--kei} is a non-empty set $X $ with maps
$$ X\to X, \, x\mapsto ax\quad {\text {and}} \quad
X\times X\to X, \, (x,y)\mapsto x*_a y\in X$$ numerated by
$a\in \alpha$ such that the following axioms are satisfied:
\begin{enumerate}
\item $ax*_a x=x$,
\item $a(x*_ay)=ax*_a ay$,
\item $(x*_a y)*_a z=(x*_a az)*_a(y*_a z),$
\item $a\tau(a)x=x$,
\item $(x*_a y)*_{\tau(a)}ay=x$,
\end{enumerate}
for all $a\in \alpha$ and $x,y,z\in X$.
Arbitrary  $\alpha$--keis can be presented by generators and relations as groups
in group theory.

\begin{Exa}
  Recall the
non-commutative ring
\begin{equation}
\Lambda=\Z[a,a_\bullet]_{a\in\alpha}/ aa_\bullet=a_\bullet a,\,a\tau(a)=1,\,a_\bullet\tau(a)_\bullet=1.
\notag
\end{equation}
Any left $\Lambda$--module $X$ becomes an $\alpha$--kei   with kei
operations   $x\mapsto ax$  and
\begin{equation} x*_a y=a_\bullet x+(1-a_\bullet
a)y. \notag
\end{equation}
The $\alpha$--keis obtained by this construction are said to be {\it abelian}.
\end{Exa}

\subsection{$\alpha$--keis of nanowords}
The theory of keis can be applied to produce homotopy invariants of
nanowords. Fix a set $\beta\subset\alpha$ such that $\tau(\beta)=\beta$.
For any nanoword $(\mathcal A,w\co\widehat {n}\to\mathcal A)$ over $\alpha$,
we define an $\alpha$--kei $\mathcal K_{\beta}(w)$. It is generated by
$n+1$ symbols $X_0,X_1,\ldots,X_n$ satisfying the following $n$ defining
relations. Each letter $A\in\mathcal A$
gives two relations. To write them down,   assume that $A$ appears in $w$
for the first time
at the $i$-th position and   for the second time at the $j$-th
position, where $i<j$. If
$a=|A| \in\beta$, then the relations are
$$
X_i =aX_{i-1},\quad
X_j =X_{j-1}*_a X_{i-1}.$$
If $a=|A|\in\alpha-\beta$, then the relations are
$$
X_i =X_{i-1}*_a X_{j-1}, \quad
X_j =aX_{j-1}.$$
The elements $V_-=X_0\in \mathcal K_{\beta}(w)$ and $V_+=X_n\in
\mathcal K_{\beta}(w)$ are called the \emph{input}
and the \emph{output}, respectively.

The idea  behind these formulas comes   from knot theory. In knot theory, every knot diagram gives rise to a
so-called quandle. Quandles are generalizations of keis and also
have only one operation, the binary operation $\ast$.
The quandle associated with a knot diagram is determined by generators,
  associated
with the arcs of the diagrams, and relations, associated with the
crossings, cf.\ the picture on the left hand side of the following
figure.\\
\begin{center}
\unitlength 0.1in
\begin{picture}( 18.3000,  6.4000)(  7.3000, -8.4000)
%
\special{pn 8}%
\special{pa 800 800}%
\special{pa 1200 410}%
\special{fp}%
\special{sh 1}%
\special{pa 1200 410}%
\special{pa 1138 442}%
\special{pa 1162 448}%
\special{pa 1166 472}%
\special{pa 1200 410}%
\special{fp}%
%
\special{pn 8}%
\special{pa 800 400}%
\special{pa 970 570}%
\special{fp}%
\special{sh 1}%
\special{pa 970 570}%
\special{pa 938 510}%
\special{pa 932 532}%
\special{pa 910 538}%
\special{pa 970 570}%
\special{fp}%
%
\special{pn 8}%
\special{pa 1040 650}%
\special{pa 1210 820}%
\special{fp}%
\special{sh 1}%
\special{pa 1210 820}%
\special{pa 1178 760}%
\special{pa 1172 782}%
\special{pa 1150 788}%
\special{pa 1210 820}%
\special{fp}%
\put(7.3000,-3.7000){\makebox(0,0)[lb]{$x$}}%
\put(11.4000,-3.8000){\makebox(0,0)[lb]{$y$}}%
\put(7.3000,-9.9000){\makebox(0,0)[lb]{$y$}}%
\put(11.5000,-10.0000){\makebox(0,0)[lb]{$x*y$}}%
%
\special{pn 8}%
\special{sh 1}%
\special{ar 1000 600 10 10 0  6.28318530717959E+0000}%
\special{sh 1}%
\special{ar 1000 600 10 10 0  6.28318530717959E+0000}%
%
%
\special{pn 8}%
\special{pa 2160 810}%
\special{pa 2560 420}%
\special{fp}%
\special{sh 1}%
\special{pa 2560 420}%
\special{pa 2498 452}%
\special{pa 2522 458}%
\special{pa 2526 482}%
\special{pa 2560 420}%
\special{fp}%
%
\special{pn 8}%
\special{pa 2160 410}%
\special{pa 2550 800}%
\special{fp}%
\special{sh 1}%
\special{pa 2550 800}%
\special{pa 2518 740}%
\special{pa 2512 762}%
\special{pa 2490 768}%
\special{pa 2550 800}%
\special{fp}%
\put(20.9000,-3.8000){\makebox(0,0)[lb]{$x$}}%
\put(25.0000,-3.9000){\makebox(0,0)[lb]{$ay$}}%
\put(20.9000,-10.0000){\makebox(0,0)[lb]{$y$}}%
\put(25.1000,-10.1000){\makebox(0,0)[lb]{$x*_a y$}}%
%
\special{pn 8}%
\special{sh 1}%
\special{ar 2360 610 10 10 0  6.28318530717959E+0000}%
\special{sh 1}%
\special{ar 2360 610 10 10 0  6.28318530717959E+0000}%
\put(23.2000,-5.3000){\makebox(0,0)[lb]{$a$}}%
\end{picture}%
\vspace{7mm}
\end{center}
In the setting of nanowords the situation is somewhat different.
First, each crossing is labeled by a letter,  $a\in \alpha$,  which
allows us to involve the operation $y\mapsto ay$ absent for knots.
The binary operation $*_{a}$ also depends on $a$. Also, the two
incoming branches are ordered. This leads us to the defining
relations as above, whose   geometric interpretation  is shown on
the right hand side of the figure.

\begin{Thm}
The     triple $(\mathcal K_{\beta}(w),
V_-,  V_+)$, considered up to isomorphism, is a homotopy
invariant of $w$.
\end{Thm}

  The $\Lambda$--module
$K_{\beta}(w)$, viewed as an $\alpha$--kei,   can be computed from
$\mathcal K_{\beta}(w)$. Namely, there is a   homomorphism of
$\alpha$--keis $\mathcal K_{\beta}(w)\to K_{\beta}(w)$ such that for
any homomorphism from $\mathcal K_{\beta}(w)$ to an abelian
$\alpha$--kei  $X$, the following diagram is commutative:
\begin{center}
\unitlength 0.1in
\begin{picture}(  9.4500,  5.7000)(  2.9000,-10.2000)
\put(2.9000,-6.2000){\makebox(0,0)[lb]{$\mathcal K_{\beta}(w)$}}%
%
\special{pn 8}%
\special{pa 700 546}%
\special{pa 1100 546}%
\special{fp}%
\special{sh 1}%
\special{pa 1100 546}%
\special{pa 1034 526}%
\special{pa 1048 546}%
\special{pa 1034 566}%
\special{pa 1100 546}%
\special{fp}%
\special{pa 1100 546}%
\special{pa 1100 546}%
\special{fp}%
\put(11.3000,-6.2000){\makebox(0,0)[lb]{$K_{\beta}(w)$}}%
\put(11.6000,-11.9000){\makebox(0,0)[lb]{$X$.}}%
%
\special{pn 8}%
\special{pa 700 646}%
\special{pa 1090 980}%
\special{fp}%
\special{sh 1}%
\special{pa 1090 980}%
\special{pa 1052 922}%
\special{pa 1050 946}%
\special{pa 1026 952}%
\special{pa 1090 980}%
\special{fp}%
%
\special{pn 8}%
\special{pa 1230 660}%
\special{pa 1230 1000}%
\special{fp}%
\special{sh 1}%
\special{pa 1230 1000}%
\special{pa 1250 934}%
\special{pa 1230 948}%
\special{pa 1210 934}%
\special{pa 1230 1000}%
\special{fp}%
\end{picture}%
\vspace{7mm}
\end{center}

\subsection{Characteristic sequences}
Consider in more detail the case   $\beta=\alpha$.  Looking at the defining relations, we easily observe that
$\mathcal K_{\beta}(w)=\mathcal K_{\alpha}(w)$ is a free $\alpha$--kei
generated by the input  $V_-$.    The output $V_+ \in \mathcal K_{\alpha}(w) $
      is a homotopy  invariant of
  $w$. The structure of  free $\alpha$--keis is poorly understood, which
  prevents us from deriving  further invariants of $w$ from $V_+$.
We focus   on a special case where more information is available.

Suppose   that the involution $\tau\co\alpha\to \alpha$ is
fixed-point-free, that is $ \tau(a)\neq a$ for all $a\in\alpha$. Fix
a   set $\alpha_+\subset \alpha$ meeting every orbit of $\tau$ in
exactly one element. Thus,
\begin{equation} \alpha=\alpha_+\cup \tau(\alpha_+),\,\,\,\alpha_+\cap \tau(\alpha_+)=\emptyset. \notag
\end{equation}
Recall the group $\Psi$  introduced in Section \ref{psiintro}. We
show   how to derive from any nanoword $w$ over $\alpha$ a finite
sequence   $(\varepsilon_1 \psi_1,\varepsilon_2
\psi_2,\ldots,\varepsilon_m \psi_m)$ with $m\geq 0$,
$\psi_1,\ldots,\psi_m\in\Psi$, and $\varepsilon_1,\ldots,
\varepsilon_m\in \{\pm 1\}$. This sequence is a homotopy invariant
of $w$ (possibly depending on $\alpha_+$). It determines
$\lambda(w)$ by
\begin{equation}
\lambda(w)=\sum_{i=1}^m\varepsilon_i\psi_i\in\Lambda=\Z\Psi.
\notag
\end{equation}
In the setting of curves, this sequence was
introduced by Silver and Williams \cite{SW2}.

We first define an $\alpha$--kei $F$ as follows. Let $F$
be the free group generated by the group $\Psi$, viewed as a set. Each  element $\psi\in \Psi$ gives rise to a  generator of   $F$, denoted
$\underline {\psi}$.  In particular,  the unit $1\in \Psi$   gives rise to a generator
$\underline{1}\in F$ which is by no means the unit of $F$.  A typical element of $F$
has the form
$$
(\underline{\psi_1})^{\varepsilon_1} (\underline{\psi_2})^{\varepsilon_2}
\cdots (\underline{\psi_m})^{\varepsilon_m}$$
where $m\geq 0$,   $
\psi_1, \psi_2,\ldots ,\psi_m\in\Psi$, and $\varepsilon_1,\ldots,\varepsilon_m\in\{\pm 1\}$. Such an element is the unit of $F$ if either $m=0$ or it can be reduced
to the case $m=0$ by applying the relations  $\underline {\psi}(\underline
{\psi})^{-1}=
(\underline {\psi})^{-1}\underline {\psi}=1$. The left action of $\Psi$ on itself extends to a group action of
$\Psi$ on $F$ by group automorphisms. The generators $a, a_\bullet \in \Psi$ act on $F$
by
\begin{equation}\begin{split}
a(&(\underline{\psi_1})^{\varepsilon_1}
(\underline{\psi_2})^{\varepsilon_2}
\cdots
(\underline{\psi_m})^{\varepsilon_m})=
(\underline{a\psi_1})^{\varepsilon_1}
(\underline{a\psi_2})^{\varepsilon_2}
\cdots
(\underline{a\psi_m})^{\varepsilon_m},\\
a_{\bullet}(&(\underline{\psi_1})^{\varepsilon_1}
(\underline{\psi_2})^{\varepsilon_2}
\cdots
(\underline{\psi_m})^{\varepsilon_m})=
(\underline{a_{\bullet}\psi_1})^{\varepsilon_1}
(\underline{a_{\bullet}\psi_2})^{\varepsilon_2}
\cdots
(\underline{a_{\bullet}\psi_m})^{\varepsilon_m}.\\
\end{split}
\notag
\end{equation}
  This defines in particular the mapping $F\to F$, $
x\mapsto ax$ for all $a\in \alpha$.
The binary operation  $x*_a y$  for  $x,y\in F$ is defined by
\begin{equation}
x*_a y=y(a_\bullet x)(a_\bullet ay)^{-1}\in F,
\notag
\end{equation}
if   $a\in\alpha_+$ and
\begin{equation}
x*_a y=(\tau(a)^{-1}_{\bullet}\tau(a)^{-1}y)^{-1}(a^{-1}_{\bullet}x)y\in F,
\notag
\end{equation}
if $a\in \alpha-\alpha_+$. These operations make $F$ into an
$\alpha$--kei.

Recall that  starting with a
nanoword $w$, we obtained a homotopy invariant  element $V_+$ of the
free $\alpha$--kei $\mathcal K_{\alpha}(w)$
on one generator $V_-$. Since $\mathcal K_{\alpha}(w)$ is free, there is a unique $\alpha$--kei
homomorphism $f:\mathcal K_{\alpha}(w)\to F$ such that $f(V_-)=\underline{1}\in F$.
Then $f(V_+)\in F$ is  a homotopy invariant of $w$.  We can expand
\begin{equation}
f(y)=(\underline{\psi_1})^{\varepsilon_1}\cdots(\underline{\psi_m})^{\varepsilon_m} \in F,
   \notag
\end{equation}
where $\psi_1,\ldots,\psi_m\in\Psi$ and $\varepsilon_1,\ldots,\varepsilon_m\in\{\pm 1\}$.
The resulting sequence
$(\varepsilon_1\psi_1,\ldots,\varepsilon_m\psi_m)$   is well-defined up to
insertion or deletion of consecutive pairs $(+\psi,-\psi)$ and $(-\psi, +\psi)$.
Deleting all such  pairs, we obtain a uniquely defined sequence
$(\varepsilon_1
\psi_1, \ldots,\varepsilon_{m'} \psi_{m'})$ with $m'\leq m$
which is a homotopy invariant of $w$.  This    is   the
\emph{characteristic sequence} of $w$.

\subsection{Examples}
1. Pick $a,b \in \alpha_+$ and consider   the  nanoword $w=ABAB$ with $ |A|=a $
and $  |B|=b  $.
It is easy to compute from the relations that  $V_+=(ba V_- *_a V_-)*_a a V_-$.
The
  characteristic sequence of $w$ is computed to be
\begin{equation}
(a,b_\bullet,b_\bullet a_\bullet ba,-b_\bullet a_\bullet
a,-b_\bullet ba). \notag
\end{equation}
In particular, if $a=b$, then this  sequence  is $(a,a_\bullet,a_\bullet^2 a^2,
-  a^2_\bullet
a,-a_\bullet aa)$.

2. Consider   the nanoword $w_1=ABACBC,\ |A|=|C|=a,\
|B|=\tau(a)\neq a$. Its   characteristic sequence (determined by any   $\alpha_+\subset \alpha$
as above such that $a\in \alpha_+$)
is:
\begin{equation}
(1,a_\bullet,-aa_\bullet,-1,a,aa_\bullet,-a^2a_\bullet,
aa_\bullet,a^2a_\bullet^2,-aa_\bullet^2,-aa_\bullet).
\notag
\end{equation}
Comparing with the previous example (for $a=b$), we obtain
that $w_1$ is not  homotopic to the nanoword $w_2=ACAC$ with $|A|=|C|=a$. This result was
claimed at the end of  Section \ref{fgi134}.

3. One might think that such a powerful invariant as the characteristic
sequence should   distinguish  arbitrary non-homotopic nanowords. However this is not
true, as shows the following example. Pick four letters $a,b,c,d\in \alpha$ (possibly
coinciding) and consider the nanoword
$$
w=ABCDCDAB,\,\,\,|A|=a,\, |B|=b,\, |C|=c,\, |D|=d.$$
An inspection shows that if $a\neq \tau (b)$ and $c\neq \tau (d)$, then the $\alpha$-pairing
of $w$ is primitive. Then $\norm {w}=4$ and $w$ is non-contractible. However, a direct
computation shows that for   $a,b\in
\alpha_+$ and
  $ c=\tau(b),  d=\tau(a)$,
the characteristic sequence of $w$
is the same as the one of  the empty nanoword.
Both consist of a single term $1\in \Psi$.

\section{Open questions and further directions}

\begin{Ques}
Classify nanowords of length $\leq 10$ up to homotopy.
\end{Ques}
In \cite{Tu1}  we give a homotopy classification of nanowords up to
length $6$. The next step is to handle the nanowords of length $8$.
Does one need new  homotopy invariants already for length $8$ ?

\begin{Ques}
Classify words of length $\leq 7$ up to homotopy.
\end{Ques}
In \cite{Tu1} we give a homotopy classification of   words up to length $5$.
One may try to classify   words by first classifying
nanowords. However, short words may desingularize into quite long nanowords.
For example, the word  $aababb$,
  desingularizes into a nanoword of length $12$.  Still,
   a classification of words of length $\leq 7$ does not
  look unrealistic because they desingularize  into a quite particular set of   nanowords.

\begin{Ques}
What (primitive) $\alpha$--pairings can be realized as
$\alpha$--pairings of nanowords?
\end{Ques}

\begin{Ques}
What   polynomials $\lambda\in\Lambda$ arise from
nanowords?
\end{Ques}
There are some simple known conditions, see \cite{Tu1}.   All new
conditions are welcome.

\begin{Ques} Is it true that all nanowords over the alphabet   consisting of a single element are
contractible ?
\end{Ques}

At the moment, nothing contradicts the conjecture that the answer is yes.

\begin{Ques}
Give a normal form for elements of a free $\alpha$--kei on one
generator.
\end{Ques}

Such a normal form (or at least  an algorithm to distinguish  elements of this $\alpha$--kei) would help to distinguish words  up to
homotopy.

One further direction is the study of cobordisms of words. Cobordism is an equivalence relation
generated by surgery on words which consists in
deleting or inserting   symmetric subwords or  subphrases. There are difficult
problems concerning  the classification of words up to cobordism. This
is studied in \cite{Tu3}.

Another interesting  direction is a study of   higher dimensional words over an alphabet
$\alpha$.  Knot theory and other
topological ideas
used above   generalize to higher dimensions.
What can be said about  similar generalizations of words?
From the topological perspective,
an $n$--dimensional nanoword is an immersion of a connected
$n$-dimensional manifold   into an $(n+1)$--dimensional manifold.
The double points of the immersion split  as a union of connected
$(n-2)$-dimensional manifolds   labeled with  letters of  $\alpha$.
The case $n=1$ is treated in the present  paper. The next
  case $n=2$ is quite mysterious. What are the appropriate analogues of the homotopy
  moves for $n=2$ ?
Although a study of high-dimensional words is tempting,
it is   hard to imagine intelligent beings communicating
with such    words.


\begin{thebibliography}{MM}




                 \bibitem[CE]{ce} G. Cairns  and D. M. Elton,
\emph{The planarity problem for signed Gauss words\/},  J. Knot Theory
Ramifications  2  (1993),    359--367.


    \bibitem[CKS]{cks} J. S. Carter, S. Kamada, M. Saito,
\emph{Stable equivalence of knots on surfaces and virtual knot cobordisms\/}.   J. Knot Theory Ramifications
11 (2002),   311--322.


                  \bibitem[CW]{cw} N. Chaves and
C. Weber,
                      \emph{Plombages de rubans et
probl\`eme des mots de
Gauss\/},     Exposition. Math.  12  (1994),
  53--77 and
124.





                  \bibitem[CR]{cr}
                    H.  Crapo and P.  Rosenstiehl,
\emph{On lacets and their manifolds\/},
Discrete Math. 233 (2001),   299--320.



                      \bibitem[DT]{dt}  C. H. Dowker
and M.B.
Thistlethwaite,
\emph{Classification of knot projections\/},
Topology Appl. 16 (1983),   19--31.





  \bibitem[Ga]{ga} C. F. Gauss,
                      Werke, Vol.\ VIII, Teubner, Leipzig, 1900, pp.
272, 282--286.


  \bibitem[GPV]{gpv} M. Goussarov, M. Polyak, and O. Viro,
\emph{Finite-type invariants of classical and virtual knots\/}.   Topology  39  (2000),
  1045--1068.

\bibitem[Kam]{Kam}   S. Kamada,
\emph{Knot invariants derived from quandles and racks\/}.  Invariants of knots and 3-manifolds (Kyoto, 2001),
103--117 (electronic), Geom. Topol. Monogr., 4, Geom. Topol. Publ., Coventry, 2002. 57M27


\bibitem[KK]{kk} N. Kamada and  S. Kamada,
\emph{Abstract link diagrams and virtual knots\/}. J. Knot Theory
Ramifications 9 (2000),  93--106.


  \bibitem[Ka]{ka} L. Kauffman,
                      \emph{Virtual knot theory\/},
                 European J. Combin.  20  (1999),    663--690.

    \bibitem[LM]{LM} L. Lov\'asz and  M. L. Marx,
\emph{A forbidden substructure characterization of Gauss codes\/}.
Acta Sci. Math. (Szeged) 38 (1976),   115--119.



              \bibitem[Ma]{Ma} M. L. Marx,
                      \emph{The Gauss realizability problem\/},
                 Proc. Amer. Math. Soc. 22 (1969), 610--613.


                      \bibitem[Ro]{ro} P. Rosenstiehl,
                      \emph{Solution alg\'ebrique du
probl\`eme de Gauss
sur la
permutation des points d'intersection d'une ou
plusieurs courbes
ferm\'ees du
plan\/},   C. R. Acad. Sci. Paris S\'er. A-B  283
(1976),
A551--A553.




   \bibitem[Sa]{SA} J. Sawollek,
  \emph{On Alexander-Conway polynomials
for virtual knots and links\/},
math.GT/9912173

\bibitem[SW1]{SW1} D. Silver and  S. Williams,
  \emph{ Polynomial invariants of virtual links\/},
J. Knot Theory Ramifications 12 (2003),   987--1000.

\bibitem[SW2]{SW2} D. Silver and  S. Williams,
  \emph{An invariant for open virtual strings\/},
   J. Knot Theory Ramifications  15  (2006),   143--152.



\bibitem[Tu1]{Tu0} V. Turaev,  \emph{Virtual strings\/},  Ann. Inst. Fourier   54
(2004),  2455--2525.
\bibitem[Tu2]{Tu1} V. Turaev, {\it Topology of words}, math.CO/0503683.
\bibitem[Tu3]{Tu2} V. Turaev, {\it Knots and words}, math.CO/math.GT/0506390.
\bibitem[Tu4]{Tu3} V. Turaev, {\it Cobordisms of words}, math.CO/0511513.
\end{thebibliography}
\end{document}